\documentclass[10pt]{article}

\usepackage{amsmath,amssymb,amsthm,amscd}

\begin{document}


\newcommand{\non}{\nonumber}
\newcommand{\scl}{\scriptstyle}
\newcommand{\wt}{\widetilde}
\newcommand{\wh}{\widehat}
\newcommand{\ot}{\otimes}
\newcommand{\fand}{\quad\text{and}\quad}
\newcommand{\Fand}{\qquad\text{and}\qquad}
\newcommand{\ts}{\,}
\newcommand{\tss}{\hspace{1pt}}
\newcommand{\tpr}{t^{\tss\prime}}
\newcommand{\spr}{s^{\tss\prime}}


\newcommand{\U}{ {\rm U}}
\newcommand{\Z}{ {\rm Z}}
\newcommand{\ZY}{ {\rm ZY}}
\newcommand{\Sr}{ {\rm S}}
\newcommand{\Ar}{ {\rm A}}
\newcommand{\Br}{ {\rm B}}
\newcommand{\Cr}{ {\rm C}}
\newcommand{\Ir}{ {\rm I}}
\newcommand{\Jr}{ {\rm J}}
\newcommand{\Prm}{ {\rm P}}
\newcommand{\X}{ {\rm X}}
\newcommand{\Y}{ {\rm Y}}
\newcommand{\SY}{ {\rm SY}}
\newcommand{\Pf}{ {\rm Pf}\ts}
\newcommand{\Hf}{ {\rm Hf}\ts}
\newcommand{\tr}{ {\rm tr}\ts}
\newcommand{\End}{{\rm{End}\ts}}
\newcommand{\Hom}{{\rm{Hom}}}
\newcommand{\sgn}{ {\rm sgn}\ts}
\newcommand{\sign}{ {\rm sign}\ts}
\newcommand{\qdet}{ {\rm qdet}\ts}
\newcommand{\sdet}{ {\rm sdet}\ts}
\newcommand{\Ber}{ {\rm Ber}\ts}
\newcommand{\inv}{ {\rm inv}\ts}
\newcommand{\inva}{ {\rm inv}}
\newcommand{\middle}{ {\rm mid} }
\newcommand{\ev}{ {\rm ev}}
\newcommand{\gr}{ {\rm gr}\ts}


\newcommand{\CC}{\mathbb{C}}
\newcommand{\ZZ}{\mathbb{Z}}


\newcommand{\Ac}{{\mathcal A}}
\newcommand{\Bc}{{\mathcal B}}
\newcommand{\Oc}{{\mathcal O}}
\newcommand{\Pc}{{\mathcal P}}
\newcommand{\Qc}{{\mathcal Q}}


\newcommand{\Sym}{\mathfrak S}
\newcommand{\h}{\mathfrak h}
\newcommand{\gl}{\mathfrak{gl}}
\newcommand{\oa}{\mathfrak{o}}
\newcommand{\spa}{\mathfrak{sp}}
\newcommand{\osp}{\mathfrak{osp}}
\newcommand{\g}{\mathfrak{g}^{}}
\newcommand{\ggot}{\mathfrak{g}}
\newcommand{\agot}{\mathfrak{a}}
\newcommand{\sll}{\mathfrak{sl}}


\newcommand{\al}{\alpha}
\newcommand{\be}{\beta}
\newcommand{\ga}{\gamma}
\newcommand{\de}{\delta^{}}
\newcommand{\ve}{\varepsilon}
\newcommand{\vp}{\varphi}
\newcommand{\la}{\lambda}
\newcommand{\La}{\Lambda}
\newcommand{\Ga}{\Gamma}
\newcommand{\si}{\sigma}
\newcommand{\ze}{\zeta}
\newcommand{\om}{\omega^{}}

\renewcommand{\theequation}{\arabic{section}.\arabic{equation}}

\newtheorem{thm}{Theorem}[section]
\newtheorem{lem}[thm]{Lemma}
\newtheorem{prop}[thm]{Proposition}
\newtheorem{cor}[thm]{Corollary}

\theoremstyle{definition}
\newtheorem{defin}[thm]{Definition}
\newtheorem{example}[thm]{Example}

\theoremstyle{remark}
\newtheorem{remark}[thm]{Remark}

\newcommand{\bth}{\begin{thm}}
\renewcommand{\eth}{\end{thm}}
\newcommand{\bpr}{\begin{prop}}
\newcommand{\epr}{\end{prop}}
\newcommand{\ble}{\begin{lem}}
\newcommand{\ele}{\end{lem}}
\newcommand{\bco}{\begin{cor}}
\newcommand{\eco}{\end{cor}}
\newcommand{\bde}{\begin{defin}}
\newcommand{\ede}{\end{defin}}
\newcommand{\bex}{\begin{example}}
\newcommand{\eex}{\end{example}}
\newcommand{\bre}{\begin{remark}}
\newcommand{\ere}{\end{remark}}

\newcommand{\bal}{\begin{aligned}}
\newcommand{\eal}{\end{aligned}}
\newcommand{\beq}{\begin{equation}}
\newcommand{\eeq}{\end{equation}}
\newcommand{\ben}{\begin{equation*}}
\newcommand{\een}{\end{equation*}}

\newcommand{\bpf}{\begin{proof}}
\newcommand{\epf}{\end{proof}}

\def\beql#1{\begin{equation}\label{#1}}

\title{\Large\bf Skew representations of twisted Yangians}

\author{{\sc A. I. Molev}\\[10mm]
School of Mathematics and Statistics\\
University of Sydney,
NSW 2006, Australia\\
{\tt alexm\hspace{0.09em}@\hspace{0.1em}maths.usyd.edu.au}
}

\date{} 


\maketitle

\vspace{15 mm}

\begin{abstract}
Analogs of the classical Sylvester theorem have been known
for matrices with entries in noncommutative algebras
including the quantized algebra of functions on
$GL_N$ and the Yangian for $\gl_N$.
We prove a version of this theorem for
the twisted Yangians $\Y(\g_N)$ associated with the orthogonal and
symplectic Lie algebras $\g_N=\oa_N$ or $\spa_N$.
This gives rise to representations of the twisted Yangian $\Y(\g_{N-M})$
on the space of homomorphisms $\Hom_{\g_M}(W,V)$,
where $W$ and $V$ are finite-dimensional irreducible
modules over $\g_M$ and $\g_N$, respectively.
In the symplectic case these representations
turn out to be irreducible and we identify them by
calculating the corresponding Drinfeld polynomials.
We also apply the quantum Sylvester theorem
to realize the twisted Yangian
as a projective limit of certain centralizers in universal
enveloping algebras.
\end{abstract}

\newpage

\section{Introduction}
\label{sec:int}
\setcounter{equation}{0}

Let $\g$ be a complex reductive Lie algebra and $\agot\subset\g$ a
reductive subalgebra.
Suppose that $V$ is a finite-dimensional irreducible $\g$-module
and consider its restriction to the subalgebra $\agot$.
This restriction
is isomorphic to a direct sum of irreducible
finite-dimensional $\agot$-modules $W_{\mu}$ with certain
multiplicities $m_{\mu}$,
\ben
V|^{}_{\agot}\cong \underset{\mu}{\oplus}\ts m_{\mu}\ts W_{\mu}.
\een
If each $W_{\mu}$ is provided with a basis and
the decomposition is multiplicity-free (i.e., $m_{\mu}\leqslant 1$ for all $\mu$)
then it can be used to get a basis of $V$ as the union of the bases of the
spaces $W_{\mu}$ which occur in the decomposition.
This observation played a key role in the
construction of the Gelfand--Tsetlin bases for the representations
of the general linear and orthogonal Lie algebras.
Although the restriction of an irreducible
finite-dimensional representation of the symplectic Lie algebra
$\spa_{2n}$ to the subalgebra $\spa_{2n-2}$
is not multiplicity-free in general,
this approach can be extended to the symplectic case
with the use of the isomorphism
\beql{dectp}
V\cong \underset{\mu}{\oplus}\ts U_{\mu}\otimes W_{\mu},
\eeq
where
\ben
U_{\mu}=\Hom_{\agot}(W_{\mu},V),\qquad \dim U_{\mu}=m_{\mu}.
\een
The space $U_{\mu}$ is an irreducible module over the
algebra $\Cr(\g,\agot)=\U(\g)^{\agot}$, the
centralizer of $\agot$
in the universal enveloping algebra $\U(\g)$;
see e.g. Dixmier~\cite[Section 9.1]{d:ae}.
Now, if some bases of the spaces $U_{\mu}$ and $W_{\mu}$ are given then
the decomposition \eqref{dectp} yields the
natural tensor product basis of $V$.
The general difficulty of this approach is the complicated
structure of the algebra $\Cr(\g,\agot)$.
For each pair of the classical Lie algebras
\ben
(\g,\agot)\quad=\quad(\gl_N,\gl_M),\quad (\oa_N,\oa_M),\quad (\spa_{N},\spa_{M}),
\een
(with even $N$ and $M$ in the symplectic case),
the centralizer $\Cr(\g,\agot)$ and its representations
can be studied with the use of
the quantum algebras called Yangians and twisted Yangians.
The {\it Yangian\/} $\Y(\gl_N)$ for the general
linear Lie algebra
$\gl_N$ is a deformation of the
universal enveloping algebra $\U(\gl_N\ot\CC[x])$
in the class of Hopf algebras; see e.g. Drinfeld~\cite{d:ha}.
The {\it twisted Yangian\/} $\Y(\g_N)$
for the orthogonal or symplectic
Lie algebra ($\g_N=\oa_N$ or $\g_N=\spa_N$)
was introduced by Olshanski~\cite{o:ty}. This is a
subalgebra of $\Y(\gl_N)$ and it can also be presented by
generators and defining relations; see also \cite{mno:yc}.
Finite-dimensional irreducible representations of the
algebras $\Y(\gl_N)$ and $\Y(\g_N)$ admit a
complete parametrization;
see Drinfeld~\cite{d:nr} and Tarasov~\cite{t:im} for the Yangian case, and
the author's work~\cite{m:fd} for the twisted Yangian case.
The Olshanski {\it centralizer construction\/} \cite{o:ri,o:ty}
provides `almost surjective' algebra homomorphisms
\beql{homcc}
\Y(\gl_{N-M})\to \Cr(\gl_N,\gl_M),\qquad \Y(\g_{N-M})\to \Cr(\g_N,\g_M)
\eeq
which allow one to equip
the corresponding $\Cr(\g,\agot)$-module $U_{\mu}$ in \eqref{dectp}
with the structure of a representation of
the Yangian or twisted Yangian, respectively.
In particular, in the case $N-M=2$ this module over the twisted Yangian
$\Y(\g_{2})$ admits a natural basis which leads to a construction of
weight bases of Gelfand--Tsetlin type for the representations
of the orthogonal and symplectic Lie algebras; see \cite{m:gtb} for
a review of these results.

In this paper we exploit the relationship between
the (twisted) Yangians and the classical Lie algebras in the reverse direction:
we use the weight bases constructed in \cite{m:gtb}
to investigate the representations of the twisted Yangians $\Y(\g_{N-M})$
emerging from the homomorphisms \eqref{homcc}.

By the results of \cite{d:nr} and \cite{m:fd}
the isomorphism class of each finite-dimensional
irreducible representation $V$ of the (twisted) Yangian is
determined by its highest weight which is a tuple of formal series
over $\CC$ in a formal parameter.
Moreover, simultaneous multiplication of all components of the highest
weight by a fixed invertible formal series corresponds to a representation
obtained from $V$ by the composition with a simple automorphism
of the (twisted) Yangian. It is natural to combine
these representations into a single {\it similarity class\/}.
In the case of the Yangian $\Y(\gl_N)$ these similarity classes
correspond to finite-dimensional
irreducible representations of the Yangian for the special linear
Lie algebra $\sll_N$. Both in the case of the Yangian and the twisted
Yangian the similarity classes are parameterized by
families of the {\it Drinfeld polynomials\/}
$(P_1(u),\dots,P_r(u))$ with some additional data in the twisted case.
Each $P_i(u)$ is a monic
polynomial in $u$, and $r$ is the rank of the corresponding Lie algebra.

Given partitions $\la=(\la_1,\dots,\la_N)$ and $\mu=(\mu_1,\dots,\mu_M)$,
let $V(\la)$ and $V(\mu)$ be the
finite-dimensional irreducible representation of $\gl_N$ and $\gl_M$
with the highest weights $\la$ and $\mu$, respectively.
The space $\Hom_{\gl_M}(V(\mu),V(\la))$ is then an irreducible
representation of the Yangian $\Y(\gl_{N-M})$. Its Drinfeld
polynomials were calculated by Nazarov and Tarasov~\cite{nt:ry}.
The result is a simple combinatorial rule which allows one to `read off'
each polynomial $P_i(u)$ from the contents of the cells of the
skew diagram $\la/\mu$. These {\it skew\/} representations of the Yangian
(they were called {\it elementary\/} in \cite{nt:ry}), may be regarded
as building blocks for the class of {\it tame\/} representations.
This class is characterized by the property that the action of
a natural commutative subalgebra of the Yangian
in such a representation is semisimple. By \cite{nt:ry},
each tame representation is isomorphic to a tensor product of
skew representations.

A different way to define the homomorphism \eqref{homcc} in the case of $\gl_N$
is provided by the {\it quantum Sylvester
theorem\/}. Recall that
the classical Sylvester theorem is the following identity
for a numerical $N\times N$ matrix $A=(a_{ij})$:
\ben
\det B=\det A\cdot \Big({a\ts}^{m+1\cdots\ts N}_{m+1\cdots\ts N}\Big)^{m-1},
\een
where $B=(b_{ij})$ is the $m\times m$ matrix formed by the minors
$
b_{ij}={a\ts}^{i,m+1\cdots\ts N}_{j,m+1\cdots\ts N}
$
of $A$. The sequences of top and bottom indices indicate the row
and column numbers of the minor, respectively.
The most general noncommutative analog of this identity
was given by Gelfand and Retakh
in the context of the theory of {\it quasideterminants\/}
originated in their work \cite{gr:dm}; see also \cite{ggrw:q} for a review
of this theory.
`Quantum' versions of this identity
apply to the matrices formed by the generators of certain quantum algebras,
and the determinants are replaced by appropriate quantum determinants.
In particular, such a version was given by Krob and Leclerc~\cite{kl:mi}
for the quantized algebra of functions on $GL_N$. Their approach is also
applicable to the Yangian $\Y(\gl_N)$. A different proof
for the Yangian case is given in \cite{m:yt} where
the corresponding quantum Sylvester theorem was used
to give a modified version of the Olshanski centralizer construction.
This provided a new definition of the skew representations of the Yangian
and the calculation of their Drinfeld polynomials.

In this paper we produce a quantum Sylvester theorem for the twisted
Yangian $\Y(\g_N)$ with the use of the Sklyanin minors
of the matrix of generators of $\Y(\g_N)$.
We first obtain the theorem for the
{\it extended twisted Yangian\/} $\X(\g_N)$
(Section~\ref{sec:ety}), following the approach of \cite{kl:mi}.
The twisted Yangian $\Y(\g_N)$ is a quotient
of $\X(\g_N)$ which yields the corresponding result for $\Y(\g_N)$
(Section~\ref{sec:ty}).
In Section~\ref{sec:sr} we apply the quantum Sylvester theorem to
construct a new homomorphism \eqref{homcc} for the twisted Yangian
and introduce the corresponding skew representations.
We show that in the symplectic case each skew representation
is irreducible and calculate its highest weight and
the Drinfeld polynomials. The Drinfeld polynomials are found by the
following simple
combinatorial rule somewhat analogous to the Yangian case \cite{nt:ry}
(see Section~\ref{sec:sr} below for a detailed formulation).
Given a partition $\nu=(\nu_1,\dots,\nu_n)$ we draw its {\it diagram\/}
$\Ga(\nu)$ as follows. First, place the row with $\nu_n$ unit cells
on the plane in such a way that the center of the leftmost cell coincides
with the origin. Then place the second row with
$\nu_{n-1}-\nu_n$ cells in such a way
that the southwest corner of this row coincides with the northeast corner
of the first row. Continuing in this manner, we complete this procedure
by placing an infinite row of cells in such a way that
its southwest corner coincides with the northeast corner
of the row with $\nu_1-\nu_2$ cells. The diagram
$\Ga(\nu)$ is obtained as the union
of the rows just placed and their images under the central symmetry
with respect to the southwest corner of the first row.
The figure below represents the diagram for the partition $\nu=(7,4)$,
where the dot indicates the origin.

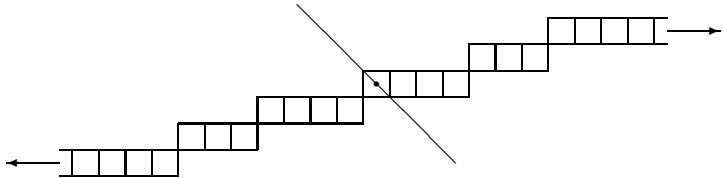
\setlength{\unitlength}{0.5em}
\begin{center}
\begin{picture}(48,14)

\put(1,0){\line(1,0){9}}
\put(1,2){\line(1,0){15}}
\put(10,4){\line(1,0){14}}
\put(16,4){\line(1,0){8}}
\put(16,6){\line(1,0){16}}
\put(24,6){\line(1,0){8}}
\put(24,8){\line(1,0){14}}
\put(32,8){\line(1,0){6}}
\put(32,10){\line(1,0){15}}
\put(38,12){\line(1,0){9}}

\multiput(2,0)(2,0){5}{\line(0,1){2}}
\multiput(10,2)(2,0){4}{\line(0,1){2}}
\multiput(16,4)(2,0){5}{\line(0,1){2}}
\multiput(24,6)(2,0){5}{\line(0,1){2}}
\multiput(32,8)(2,0){4}{\line(0,1){2}}
\multiput(38,10)(2,0){5}{\line(0,1){2}}

\put(1,1){\vector(-1,0){4}}
\put(47,11){\vector(1,0){4}}

\put(25,7){\circle*{0.3}}

\put(19,13){\line(1,-1){12}}

\end{picture}
\end{center}
\setlength{\unitlength}{1pt}

\bigskip
\noindent
To each cell of the diagram
we attach its diagonal number, where by diagonals we mean
the lines passing northwest-southeast through the integer points
of the plane.
The line on the figure indicates the $0$-th diagonal and the diagonal numbers
are consecutive integers increasing from right to left.
For any nonnegative integer $p$
denote by $\Ga(\la)^{(p)}$
the diagram $\Ga(\la)$ lifted $p$ units up.

Suppose now that $V(\la)$ and $V(\mu)$ are the irreducible
finite-dimensional representations
of $\spa_{2n}$ and $\spa_{2m}$ corresponding to partitions
$\la$ and $\mu$ having $n$ and $m$ parts, respectively.
Then the Drinfeld
polynomials $P_1(u),\dots, P_{n-m}(u)$ for the
skew representation $\Hom_{\spa_{2m}}(V(\mu),V(\la))$
of the twisted Yangian
$\Y(\spa_{2n-2m})$ can be calculated by the following rule:
all roots of the polynomial $P_k(u)$ are simple and they
coincide with the diagonal numbers decreased by $1/2$ of the cells of the
intersection $\Ga(\mu)\cap \Ga(\la)^{(k-1)}$
(see Theorem~\ref{thm:dp} and Example~\ref{ex:lamu} below).

Finally, in Section~\ref{sec:cc} we give a realization of the twisted Yangian
$\Y(\g_N)$ as a projective limit of centralizers
in the universal enveloping algebras.
This is a new version of the centralizer construction (cf. \cite{o:ty,mo:cc})
which is based on the quantum Sylvester theorem.

The recent work of Nazarov~\cite{n:rt} is also devoted
to the skew representations of the twisted Yangians although from
a different perspective.
He uses the classical Weyl's approach
and gives a realization of the skew representations in the tensor powers
of the vector representation by applying certain generalized Young symmetrizers.

\section{Extended twisted Yangian}
\label{sec:ety}
\setcounter{equation}{0}

We start by stating and proving some auxiliary results
about the extended twisted Yangian
$\X(\g_N)$; see \cite{mno:yc} for more details.

\subsection{Preliminaries}

We shall be considering the orthogonal and symplectic cases
simultaneously, unless otherwise stated.
Given a positive integer $N$, we number the rows and columns
of $N\times N$ matrices by the indices
$\{-n,\dots,-1,0,1,\dots, n\}$ if $N=2n+1$, and by
$\{-n,\dots,-1,1,\dots, n\}$ if $N=2n$.
Similarly, in the latter case the range of
summation indices $-n\leqslant i,j\leqslant n$
will usually exclude $0$.
It will be convenient to use
the symbol $\theta_{ij}$ which is defined by
\ben
\theta_{ij}=\begin{cases} 1\quad&\text{in the orthogonal case},\\
\sgn i\cdot \sgn j\quad&\text{in the symplectic case}.\end{cases}
\een
Throughout the paper, whenever the double sign $\pm{}$ or $\mp{}$ occurs,
the upper sign corresponds to the orthogonal case and the lower sign to
the symplectic case.
By $A\mapsto A^t$ we will denote the matrix
transposition such that
$
(A^t)_{ij}=\theta_{ij}\ts A_{-j,-i}.
$
Let the $E_{ij}$ denote the standard basis vectors of the
general linear Lie algebra $\gl_N$.
These vectors may be also regarded as elements of the universal enveloping
algebra $\U(\gl_N)$. For this reason
we want to distinguish the $E_{ij}$ from
the standard matrix units $e_{ij}$
which are considered as basis elements
of the endomorphism algebra $\End\CC^N$.
Introduce the following elements of the Lie algebra $\gl_N$:
\ben
F_{ij}=E_{ij}-\theta_{ij}E_{-j,-i},\qquad -n\leqslant i,j\leqslant n.
\een
The Lie subalgebra $\g_N$ of $\gl_N$ spanned
by the elements $F_{ij}$ is isomorphic to
the orthogonal Lie algebra $\oa_N$ or the symplectic Lie algebra
$\spa_N$ (in the latter case $N$ is even).

The {\it extended twisted Yangian\/} $\X(\g_N)$
corresponding to the Lie algebra $\g_N$ is the associative algebra
with generators $s_{ij}^{(1)},\ s_{ij}^{(2)},\dots$ where
$-n\leqslant i,j\leqslant n$,
subject to the defining relations written in
terms of the generating series
\ben
s_{ij}(u)=\delta_{ij}+ s_{ij}^{(1)}u^{-1}+s_{ij}^{(2)}u^{-2}+\cdots
\in\X(\g_N)[[u^{-1}]]
\een
as follows
\ben
\bal
(u^2-v^2)\ts
[s_{ij}(u),s_{kl}(v)]=(u+v)\ts&\big(s_{kj}(u)s_{il}(v)-s_{kj}(v)s_{il}(u)\big)\\
{}-(u-v)\ts&\big(\theta_{k,-j}s_{i,-k}(u)s_{-j,l}(v)-
\theta_{i,-l}s_{k,-i}(v)s_{-l,j}(u)\big)\\
{}+\theta_{i,-j}&\big(s_{k,-i}(u)s_{-j,l}(v)-s_{k,-i}(v)s_{-j,l}(u)\big),
\eal
\een
where $u$ and $v$ denote formal variables.
The defining relations can also be presented in a convenient
matrix form.
Denote by $S(u)$ the $N\times N$ matrix whose
$ij$-th entry is $s_{ij}(u)$. We may regard
$S(u)$ as an element of the algebra $\X(\g_N)[[u^{-1}]]\ot\End\CC^N$ given by
\ben
S(u)=\sum_{i,j} s_{ij}(u)\ot e_{ij},
\een
where the $e_{ij}$ denote the standard matrix
units. For any positive integer $m$ we shall be using the algebras of the form
\beql{multitp}
\X(\g_N)[[u^{-1}]]\ot \End \CC^N\ot \cdots\ot \End \CC^N,
\eeq
with $m$ copies of $\End \CC^N$. For any $a\in\{1,\dots,m\}$ we denote by $S_a(u)$
the matrix $S(u)$ which acts on the $a$-th copy of $\End \CC^N$. That is,
$S_a(u)$ is an element of the algebra \eqref{multitp} of the form
\ben
S_a(u)=\sum_{i,j} s_{ij}(u)\ot 1\ot \cdots \ot 1\ot e_{ij}\ot 1\ot \cdots\ot 1,
\een
where the $e_{ij}$ belong to the $a$-th copy of $\End \CC^N$ and
$1$ is the identity matrix. Similarly, if
\ben
C=\sum_{i,j,k,l} c^{}_{ijkl}\ts e_{ij}\ot e_{kl}\in \End \CC^N\ot\End \CC^N,
\een
then for distinct indices $a,b\in\{1,\dots,m\}$ we introduce
the element $C_{ab}$ of the algebra \eqref{multitp} by
\ben
C_{ab}=\sum_{i,j,k,l} c^{}_{ijkl}\ts
1\ot 1\ot \cdots \ot 1\ot e_{ij}\ot 1\ot \cdots\ot 1\ot e_{kl}\ot 1\ot \cdots\ot 1,
\een
where the  $e_{ij}$ and $e_{kl}$ belong to the $a$-th and $b$-th copies of
$\End \CC^N$, respectively.
Consider now the permutation operator
\ben
P=\sum_{i,j} e_{ij}\ot e_{ji}\in \End \CC^N\ot\End \CC^N.
\een
The rational function
$
R(u)=1-Pu^{-1}
$
with values in the tensor product algebra $\End \CC^N\ot\End \CC^N$ is called
the {\it Yang\/} $R$-{\it matrix\/}.
Introduce its transposed $R^{\tss t}(u)$ by
\beql{rt}
R^{\tss t}(u)=1-Q\ts u^{-1}, \qquad Q
=\sum_{i,j}\theta_{ij}\ts e_{-j,-i}\ot e_{ji}.
\eeq
The defining relations for the extended
twisted Yangian $\X(\g_N)$ are equivalent to
the {\it quaternary relation\/}
\beql{quater}
R(u-v)S_1(u)R^{\tss t}(-u-v)S_2(v)=S_2(v)R^{\tss t}(-u-v)S_1(u)R(u-v).
\eeq

Let $u_1,\dots,u_k$ be independent variables.
For $k\geqslant 2$ consider the rational function
$R(u_1,\dots,u_k)$ with values in
$(\End\CC^N)^{\ot k}$ defined by
\ben
R(u_1,\dots,u_k)=(R_{k-1,k})(R_{k-2,k}R_{k-2,k-1}) \cdots (R_{1k}
\cdots
R_{12}),
\een
where
we abbreviate $R_{ij}=R_{ij}(u_i-u_j)$.
Set
\ben
S_i=S_i(u_i),\quad 1\leqslant i\leqslant k \qquad \text{and}\qquad
R^{\tss t}_{ij}=R^{\tss t}_{ji}=R^{\tss t}_{ij}(-u_i-u_j),
\quad 1\leqslant i<j \leqslant k.
\een
For an arbitrary permutation $(p_1,\ldots,p_k)$ of the numbers
$1,\ldots,k$, we abbreviate
\ben
\langle S_{p_1},\ldots,S_{p_k}\rangle=S_{p_1}(R^{\tss t}_{p_1p_2}\cdots
R^{\tss t}_{p_1p_k}) S_{p_2}(R^{\tss t}_{p_2p_3}\cdots
R^{\tss t}_{p_2p_k})\cdots S_{p_k}.
\een
The identity
\beql{fundam}
R(u_1,\ldots,u_k)\langle S_1,\ldots,S_k\rangle=\langle
S_k,\ldots,S_1\rangle R(u_1,\ldots,u_k)
\eeq
can be deduced from the quaternary relation \eqref{quater}; see
\cite[Proposition~4.2]{mno:yc}.
Now specialize the variables $u_i$ by setting
\beql{speci}
u_i=u-i+1,\qquad i=1,\dots,k.
\eeq
It is well known that
under this specialization, $R(u_1,\ldots,u_k)$ coincides with the
anti-symmetrization operator $A_k$ on $(\CC^N)^{\ot k}$, where
\ben
A_k=\sum_{\sigma\in\Sym_k}\sgn\sigma\cdot P_{\sigma},
\een
and $P_{\sigma}$ denotes the image of $\sigma\in\Sym_k$ under the natural
action of $\Sym_k$ on $(\CC^N)^{\ot k}$; see e.g. \cite[Proposition~2.3]{mno:yc}.
Hence specializing the variables
in \eqref{fundam} we get
\ben
A_k \ts\langle S_1,\ldots,S_k\rangle=\langle
S_k,\ldots,S_1\rangle\ts A_k.
\een
This element of the tensor product
$\X(\g_N)[[u^{-1}]]\ot(\End\CC^N)^{\ot k}$ can be written as
\ben
\sum {s\ts}^{a_1\cdots\ts a_k}_{b_1\cdots\ts b_k}(u)\ot e_{a_1b_1}\ot \cdots
\ot e_{a_kb_k},
\een
summed over the indices $a_i,b_i\in\{-n,\dots,n\}$.
We also set
${s\tss}^{a}_{b}(u)=s_{ab}(u)$.
We call the elements
${s\ts}^{a_1\cdots\ts a_k}_{b_1\cdots\ts b_k}(u)$ of $\X(\g_N)[[u^{-1}]]$
the {\it Sklyanin minors\/} of the matrix $S(u)$.
Clearly, the Sklyanin minors are skew-symmetric
with respect to permutations of the upper indices and of the lower indices:
\ben
{s\ts}^{a_{\si(1)}\cdots\ts a_{\si(k)}}_{b_1\cdots\ts b_k}(u)
=\sgn \si\cdot{s\ts}^{a_1\cdots\ts a_k}_{b_1\cdots\ts b_k}(u)
\quad\text{and}\quad
{s\ts}^{a_1\cdots\ts a_k}_{b_{\si(1)}\cdots\ts b_{\si(k)}}(u)
=\sgn \si\cdot{s\ts}^{a_1\cdots\ts a_k}_{b_1\cdots\ts b_k}(u)
\een
for any $\si\in\Sym_k$.

\bpr\label{prop:comsm}
We have the relations
\ben
\bal
&(u^2-v^2)\ts
[s_{pq}(u),{s\ts}^{a_1\cdots\ts a_k}_{b_1\cdots\ts b_k}(v)]\\
&{}={}(u+v)\ts\sum_{i=1}^k\Big(
{s}^{}_{a_iq}(u)\ts
{s\ts}^{a_1\cdots\ts p\ts\cdots\ts a_k}_{b_1\ \cdots\ \  b_k}(v)
-{s\ts}^{a_1\ \cdots\ \  a_k}_{b_1\cdots\ts
q\ts\cdots\ts b_k}(v)\ts{s}^{}_{pb_i}(u)\Big)\\
&{}-{}(u-v)\ts\sum_{i=1}^k\Big(\theta_{a_i,-q}\ts
{s}^{}_{p,-a_i}(u)\ts
{s\ts}^{a_1\cdots\ts {-q}\ts\cdots\ts a_k}_{b_1\ \cdots\ \  b_k}(v)
-\theta_{p,-b_i}\ts{s\ts}^{a_1\ \cdots\ \  a_k}_{b_1\cdots\ts
{-p}\ts\cdots\ts b_k}(v)\ts{s}^{}_{-b_i,q}(u)\Big)\\
&{}+{}\theta_{p,-q}\ts\sum_{i=1}^k\Big(
{s}^{}_{a_i,-p}(u)\ts
{s\ts}^{a_1\cdots\ts {-q}\ts\cdots\ts a_k}_{b_1\ \cdots\ \  b_k}(v)
-{s\ts}^{a_1\ \cdots\ \  a_k}_{b_1\cdots\ts
{-p}\ts\cdots\ts b_k}(v)\ts{s}^{}_{-q,b_i}(u)\Big)\\
&{}+{}\sum_{i\ne j}\Big(\theta_{a_j,-q}\ts
{s}^{}_{a_i,-a_j}(u)\ts
{s\ts}^{a_1\cdots\ts p\ts\cdots\ts {-q}\ts\cdots\ts a_k}_{b_1\ \cdots\ \  b_k}(v)
-\theta_{p,-b_i}{s\ts}^{a_1\ \cdots\ \  a_k}_{b_1\cdots\ts
{-p}\ts\cdots\ts{q}\ts\cdots\ts b_k}(v)\ts{s}^{}_{-b_i,b_j}(u)\Big),
\eal
\een
where in the Sklyanin minors the indices $p$ and $q$
replace $a_i$ and $b_i$, respectively, in the first sum;
the indices $-q$ and $-p$
replace $a_i$ and $b_i$, respectively, in the second and third sums;
in the fourth sum $p$ and $-q$ replace $a_i$ and $a_j$, respectively,
and $-p$ and $q$ replace $b_i$ and $b_j$, respectively.
\epr

\bpf
By \eqref{fundam}, we have the relation
\begin{multline}\label{twruvv}
R(u,v,v-1,\dots, v-k+1)\ts \langle S_0,\ldots,S_k\rangle\\
=\langle S_k,\ldots,S_0\rangle\ts R(u,v,v-1,\dots, v-m+1),
\end{multline}
where we have used an extra copy of the algebra $\End\CC^N$ labelled
by $0$ and the parameters are specialized as follows
\ben
u_0=u, \fand u_i=v-i+1\quad\text{for}\quad i=1,\dots,k.
\een
Then one easily verifies (see e.g. \cite{m:yt}) that the product
of $R$-matrices in \eqref{twruvv} simplifies to
\ben
R(u,v,v-1,\dots, v-k+1)=A_k\ts \Big(1-\frac{1}{u-v}(P_{01}+\cdots+P_{0k})\Big).
\een
Applying the transposition over the zeroth copy of $\End\CC^N$
and replacing $u$ by $-u$ we also deduce
that
\beql{araq}
A_k\ts R^{\tss t}_{01}\cdots R^{\tss t}_{0k}
 =A_k\ts \Big(1+\frac{1}{u+v}(Q_{01}+\cdots+Q_{0k})\Big).
\eeq
Hence \eqref{twruvv} takes the form
\begin{multline}
\Big(1-\frac{1}{u-v}(P_{01}+\cdots+P_{0k})\Big)\ts S_0(u)
\ts \Big(1+\frac{1}{u+v}(Q_{01}+\cdots+Q_{0k})\Big)
\ts A_k\ts \langle S_1,\ldots,S_k\rangle\\
=\langle S_k,\ldots,S_1\rangle\ts A_k\ts
\Big(1+\frac{1}{u+v}(Q_{01}+\cdots+Q_{0k})\Big)\ts S_0(u)
\ts \Big(1-\frac{1}{u-v}(P_{01}+\cdots+P_{0k})\Big).
\non
\end{multline}
It remains to apply both sides to the vector $e_q\ot e_{b_1}\ot\cdots\ot e_{b_k}$
and compare the coefficients at the vector $e_p\ot e_{a_1}\ot\cdots\ot e_{a_k}$,
where the $e_i$ denote the canonical basis vectors of $\CC^N$.
\epf

\bco\label{cor:centsm}
Suppose that for some indices $i,j,l,m\in\{1,\dots,k\}$ we have
$a_i=-b_l$ and $b_j=-a_m$. Then
\ben
[s_{a_ib_j}(u),{s\ts}^{a_1\cdots\ts a_k}_{b_1\cdots\ts b_k}(v)]=0.
\een
\eco

\bpf
By the skew-symmetry property, the Sklyanin minor is zero if it
has two repeated upper or lower indices. Hence we may assume that
$i=m$ if and only if $j=l$. Suppose first that $i\ne m$. Then
using the skew-symmetry of Sklyanin minors, we derive from
Proposition~\ref{prop:comsm} that
\ben
(u-v-1)(u+v+1)\ts
[s_{a_ib_j}(u),{s\ts}^{a_1\cdots\ts a_k}_{b_1\cdots\ts b_k}(v)]
=\theta_{a_i,-b_j}\ts
[s_{-b_j,-a_i}(u),{s\ts}^{a_1\cdots\ts a_k}_{b_1\cdots\ts b_k}(v)].
\een
The same relation holds with $i$ and $j$ replaced by $m$ and $l$,
respectively, which proves the claim in the case under consideration.
If $a_i=-b_j$ then Proposition~\ref{prop:comsm} immediately gives
\ben
(u-v-1)(u+v+1)\ts
[s_{a_ib_j}(u),{s\ts}^{a_1\cdots\ts a_k}_{b_1\cdots\ts b_k}(v)]
=0,
\een
completing the proof.
\epf

The series
\ben
\sdet S(u)={s}^{-n\ts\cdots\ts n}_{-n\ts\cdots\ts n}(u)\in \X(\g_N)[[u^{-1}]]
\een
is called the {\it Sklyanin determinant\/} of the matrix $S(u)$.
Corollary~\ref{cor:centsm} implies that
all the coefficients of this series belong to the center
of the algebra $\X(\g_N)$; see also \cite[Theorem~4.8]{mno:yc}
for a slightly different proof.

The matrix $S(u)$ is invertible and we shall denote by $S^{-1}(u)$
the inverse matrix.
The mapping
\beql{invauto}
\varpi_N:S(u)\mapsto S^{-1}\big({-u}-{N}/{2}\big)
\eeq
defines
an involutive automorphism of the algebra $\X(\g_N)$;
see \cite[Proposition~6.5]{mno:yc}.

The {\it Sklyanin comatrix\/} $\wh S(u)$ is defined by
the relation
\beql{comdef}
\wh S(u)\ts S(u-N+1)=\sdet S(u).
\eeq
Due to \eqref{invauto}, the mapping
\beql{comauto}
S(u)\mapsto \wh S(-u+N/2-1)
\eeq
defines a homomorphism of $\X(\g_N)$ into itself.

We shall also use the {\it auxiliary minors\/}
${\check s\ts}^{a_1\cdots\ts a_k}_{b_1\cdots\ts b_{k-1},c}(u)
\in \X(\g_N)[[u^{-1}]]$
defined by
\begin{multline}\label{auxmin}
A_k \ts\langle S_1,\ldots,S_{k-1}\rangle\ts R^{\tss t}_{1k}
\cdots R^{\tss t}_{k-1,k}\\
=\sum {\check s\ts}^{a_1\cdots\ts a_k}_{b_1\cdots\ts b_{k-1},c}(u)
\ot e_{a_1b_1}\ot \cdots
\ot e_{a_{k-1}b_{k-1}}\ot e_{a_{k}c},
\end{multline}
summed over $a_i,b_i,c\in\{-n,\dots,n\}$.
Since
\ben
A_k \ts\langle S_1,\ldots,S_{k-1}\rangle\ts R^{\tss t}_{1k}
\cdots R^{\tss t}_{k-1,k}\ts S_k=A_k \ts\langle S_1,\ldots,S_{k}\rangle,
\een
we immediately obtain the relation
\beql{sscheck}
{s\ts}^{a_1\cdots\ts a_k}_{b_1\cdots\ts b_k}(u)
=\sum_{c=-n}^n
{\check s\ts}^{a_1\cdots\ts a_k}_{b_1\cdots\ts b_{k-1},c}(u)
\ts s_{c\tss b_k}(u-k+1).
\eeq
We obviously have
\ben
{\check s\ts}^{a_{\si(1)}\cdots\ts a_{\si(k)}}_{b_1\cdots\ts b_{k-1},c}(u)
=\sgn \si\cdot{\check s\ts}^{a_1\cdots\ts a_k}_{b_1\cdots\ts b_{k-1},c}(u)
\een
for any $\si\in\Sym_k$. Also,
\ben
{\check s\ts}^{a_1\cdots\ts a_k}_{b_{\si(1)}\cdots\ts b_{\si(k-1)},c}(u)
=\sgn \si\cdot{\check s\ts}^{a_1\cdots\ts a_k}_{b_1\cdots\ts b_{k-1},c}(u)
\een
for any $\si\in\Sym_{k-1}$; see \cite{m:sd}.
Furthermore, it is
straightforward to obtain the following property of the auxiliary
minors from their definition (cf. \cite[Proposition~4.4]{m:sd}):
if $c\notin\{a_1,\dots,a_{k-1}\}$
and $c\notin\{-b_1,\dots,-b_{k-1}\}$ then
\beql{auzero}
{\check s\ts}^{a_1\cdots\ts a_k}_{b_1\cdots\ts b_{k-1},c}(u)=0
\eeq
if $c\ne a_k$, while
\beql{aucc}
{\check s\ts}^{a_1\cdots\ts a_{k-1},\ts c}_{b_1\cdots\ts b_{k-1},\ts c}(u)=
{s\ts}^{a_1\cdots\ts a_{k-1}}_{b_1\cdots\ts b_{k-1}}(u).
\eeq
Set $(a_1,\dots,a_N)=(-n,\dots,n)$. Then
the matrix elements $\wh s_{a_ia_j}(u)$ of the Sklyanin
comatrix $\wh S(u)$ are given by
\beql{comam}
\wh s_{a_ia_j}(u)=(-1)^{N-i}\ts
{\check s\ts}^{a_1\cdots\ts a_N}_{a_1\cdots\tss
\wh a_i\tss\cdots \ts a_N,\ts a_j}(u),
\eeq
where the hat on the right hand side indicates the index to be omitted;
see \cite[Section~6]{m:sd}.

\subsection{Sylvester theorem}

We shall need
some complementary minor identities for the algebra $\X(\g_N)$;
cf. \cite{bk:pp}, \cite{kl:mi}.
Fix a nonnegative integer $m<n$ and set $M=2m$ or $M=2m+1$
if $N=2n$ or $N=2n+1$, respectively, so that $N-M=2n-2m$.
Set $\Bc=\{-m,\dots,m\}$ and
denote by $\Ac$ the complement of the subset $\Bc$ in the set
$\{-n,\dots,n\}$ so that the elements of the set $\Ac$ are
\ben
(a^{}_1,\dots,a^{}_{N-M})=(-n,\dots,-m-1,m+1,\dots,n).
\een
The images of some Sklyanin minors and auxiliary minors
with respect to the automorphism $\varpi_N$ defined in \eqref{invauto}
are provided by the following.

\bpr\label{prop:blockstw}
We have the identities
\beql{blockstw}
\sdet S(u)\cdot
\varpi_N\Big({s\ts}^{a_1\cdots\ts a_{N-M}}_{a_1\cdots\ts a_{N-M}}
(-u+N/2-1)\Big)
={s}^{-m\ts\cdots\ts m}_{-m\ts\cdots\ts m}(u),
\eeq
and
\begin{multline}\label{blockstwsv}
\sdet S(u)\cdot
\varpi_N\Big({\check s\ts}^{a_1\cdots\ts a_{N-M}}_{a_1\cdots\tss
\wh a_i\tss\cdots \ts a_{N-M},\ts a_j}(-u+N/2-1)\Big)\\
=(-1)^{N-M-i}\ts
{s}^{-m\ts\cdots\ts m,\ts a_i}_{-m\ts\cdots\ts m,\ts a_j}(u),
\end{multline}
for any $i,j\in\{1,\dots,N-M\}$. In particular,
\beql{sdetcirc}
\sdet S(u)\cdot\varpi_N\big(\sdet S(-u+N/2-1)\big)=1.
\eeq
\epr

\bpf
By definition of the Sklyanin determinant,
\beql{bydefsdet}
A_N \ts\langle S_1,\ldots,S_N\rangle=A_N\ts\sdet S(u).
\eeq
This implies the relation
\begin{multline}
A_N \ts\langle S_1,\ldots,S_{M}\rangle\prod_{i=1,\dots,M}^{\rightarrow}
(R^{\tss t}_{i,M+1}\cdots R^{\tss t}_{iN})\\
{}=A_N\ts\sdet S(u)\ts S^{-1}_N(R^{\tss t}_{N-1,N})^{-1}
\ts S^{-1}_{N-1}\cdots \ts (R^{\tss t}_{M+1,M+2})^{-1} S^{-1}_{M+1}.
\non
\end{multline}
Note that since $Q^2=NQ$ we find from \eqref{rt} that
\ben
R^{\tss t}(u)^{-1}=R^{\tss t}(-u+N).
\een
Therefore, the right hand side can be written as
\ben
A_N\ts\sdet S(u)\ts S^{\tss\circ}_N
R^{\tss\circ}_{N-1,N}
\ts S^{\tss\circ}_{N-1}\cdots \ts R^{\tss\circ}_{M+1,M+2}\ts S^{\tss\circ}_{M+1},
\een
where we have used the notation $S^{\tss\circ}(u)=\varpi_N(S(u))$ and
\beql{ucircnew}
S^{\tss\circ}_i=S^{\tss\circ}_i(u^{\tss\circ}_i),\qquad R^{\tss\circ}_{ij}=
R^{\tss t}_{ij}(-u^{\tss\circ}_i-u^{\tss\circ}_j),\qquad u^{\tss\circ}_i=-u_i-N/2
\eeq
with $u_i=u-i+1$ for $i=1,\dots,N$. Now apply both sides to the basis vector
\ben
v=e_{-m}\ot\cdots\ot e_m\otimes e_{a_1}\ot
\cdots\ot e_{a_{N-M}},
\een
where the $e_r$ denote the canonical basis vectors of $\CC^N$.
Comparing the coefficients
at the vector $A_N\ts v$ we come to \eqref{blockstw}.

The proof of \eqref{blockstwsv} is similar with the use of the relation
\begin{multline}
\label{asinvnew}
A_N \ts\langle S_1,\ldots,S_{M+1}\rangle\prod_{i=1,\dots,M}^{\rightarrow}
(R^{\tss t}_{i,M+2}\cdots R^{\tss t}_{iN})\\
{}=A_N\ts\sdet S(u)\ts S^{-1}_N(R^{\tss t}_{N-1,N})^{-1}
\ts S^{-1}_{N-1}\cdots \ts S^{-1}_{M+2}(R^{\tss t}_{M+1,N})^{-1}
\cdots(R^{\tss t}_{M+1,M+2})^{-1}
\end{multline}
implied by \eqref{bydefsdet}. Using \eqref{ucircnew} we can
rewrite the right hand side as
\ben
A_N\ts\sdet S(u)\ts S^{\tss\circ}_N
R^{\tss\circ}_{N-1,N}
\ts S^{\tss\circ}_{N-1}\cdots \ts S^{\tss\circ}_{M+2}R^{\tss\circ}_{M+1,N}
\cdots R^{\tss\circ}_{M+1,M+2}.
\een
The proof is completed by applying both sides to the vector
\ben
v^{}_{ij}=e_{-m}\ot\cdots\ot e_m\otimes e_{a_j}\otimes e_{a_1}\ot
\cdots\otimes \wh e_{a_i}\otimes\cdots
\otimes e_{a_{N-M}}
\een
and comparing the coefficients
at the vector $A_N\ts v^{}_{11}$.
\epf

It follows from the defining relations for the extended twisted
Yangian that
the subalgebra of $\X(\g_N)$ generated by the elements $s_{ij}^{(r)}$
with $i,j\in\Bc$ can be regarded as a natural homomorphic image of
the extended twisted Yangian $\X(\g_{M})$. The homomorphism
takes the generators $s_{ij}^{(r)}$ of $\X(\g_{M})$
to the elements of $\X(\g_N)$ with the same name.
We let $\X(\g_{N-M})$ denote the
extended twisted Yangian whose generator series $s_{ab}(u)$
are enumerated by elements $a,b\in\Ac$. The mapping
$
\upsilon:\X(\g_{N-M})\mapsto \X(\g_N)
$
which sends $s_{ab}(u)$ to the series with the same name in
$\X(\g_N)$ is an algebra homomorphism.
Following \cite{bk:pp} consider the
homomorphism $\vartheta:\X(\g_{N-M})\mapsto \X(\g_N)$
defined as the composition
\beql{vartheta}
\vartheta=\varpi_{N}\circ \upsilon\circ \varpi_{N-M}.
\eeq
Its action on the generators of $\X(\g_{N-M})$ can be described
with the use of the Gelfand--Retakh quasideterminants \cite{gr:dm, gr:tn}.
Given an arbitrary
$N\times N$ matrix $X=(x_{ij})$ over a ring with $1$,
denote by $X^{ij}$ the matrix obtained from $X$ by deleting the $i$-th row
and $j$-th column. Suppose that the matrix
$X^{ij}$ is invertible.
The $ij$-{\em th quasideterminant
of} $X$ is defined by the formula
\ben
|X|_{ij}=x_{ij}-\sum_{k\ne j,l\ne i}
x_{ik}\big((X^{ij})^{-1}\big)_{kl}\ts x_{lj}.
\een
In a more graphic fashion, the quasideterminant $|X|_{ij}$ is denoted
by boxing the entry $x_{ij}$. The following proposition
is proved in the same way as its Yangian counterpart;
see \cite[Lemma~4.2]{bk:pp}.

\bpr\label{prop:psiactstw}
For any $1\leqslant i,j\leqslant N-M$, we have
\ben
\vartheta: s_{a_ia_j}(u+M/2)\mapsto
\left|\begin{matrix} s_{-m,-m}(u)&\cdots&s_{-m,m}(u)&s_{-m,a_j}(u)\\
                          \vdots&\ddots&\vdots&\vdots\\
                         s_{m,-m}(u)&\cdots&s_{m,m}(u)&s_{m,a_j}(u)\\
               s_{a_i,-m}(u)&\cdots&s_{a_i,m}(u)&\boxed{s_{a_i,a_j}(u)}
\end{matrix}\right|.
\een
\epr

We shall need an expression for the image of $s_{a_ia_j}(u)$
in terms of Sklyanin minors.

\bpr\label{prop:psimskmin}
For any $1\leqslant i,j\leqslant N-M$, we have
\ben
\vartheta: s_{a_ia_j}(u)\mapsto
\Big[{s\ts}^{-m\dots\ts m}_{-m\dots\ts m}(u+M/2)\Big]^{-1}\cdot
{s\ts}^{-m\dots\ts m,\ts a_i}_{-m\dots\ts m,\ts a_j}(u+M/2).
\een
\epr

\bpf
Considering \eqref{blockstwsv} as a relation
in the algebra $\X(\g_{N-M})$ by replacing $N$ with $N-M$ and
replacing $M$ with $0$,
we get
\begin{multline}
{s\ts}^{a_1\cdots\ts a_{N-M}}_{a_1\cdots\ts a_{N-M}}(u)\cdot
\varpi_{N-M}\Big({\check s\ts}^{a_1\cdots\ts a_{N-M}}_{a_1\cdots\tss
\wh a_i\tss\cdots \ts a_{N-M},\ts a_j}(-u+n-m-1)\Big)
\non\\
=(-1)^{N-M-i}\ts
{s}_{a_ia_j}(u).
\end{multline}
Now apply $\varpi_{N-M}$ to both sides
and use \eqref{sdetcirc} to obtain
\begin{multline}
\varpi_{N-M}:{s}_{a_ia_j}(u)\mapsto (-1)^{N-M-i}
\Big[{s\ts}^{a_1\cdots\ts a_{N-M}}_{a_1\cdots\ts a_{N-M}}
(-u+n-m-1)\Big]^{-1}
\non\\
{}\times {\check s\ts}^{a_1\cdots\ts a_{N-M}}_{a_1\cdots\tss
\wh a_i\tss\cdots \ts a_{N-M},\ts a_j}(-u+n-m-1).
\end{multline}
Next apply the homomorphism $\upsilon$ and observe that
the images of the Sklyanin minors and auxiliary minors
occurring in this expression coincide with
the respective minors in the algebra $\X(\g_N)$.
Indeed, we have the following
identity analogous to \eqref{araq} which is verified in the same way:
\beql{artrtsi}
A_{k-1} \ts R^{\tss t}_{1k}
\cdots R^{\tss t}_{k-1,k}=A_{k-1}\Big(1+\frac{Q_{1k}+\cdots+Q_{k-1,k}}{2u-k+1}\Big).
\eeq
Therefore the left hand side of \eqref{auxmin}
can be written as
\ben
A_k\ts S_1\ts \Big(1+\frac{Q_{12}}{2u-1}\Big)\ts S_2
\ts \Big(1+\frac{Q_{13}+Q_{23}}{2u-2}\Big)\ts\cdots S_{k-1}\ts
\Big(1+\frac{Q_{1k}+\cdots+Q_{k-1,k}}{2u-k+1}\Big).
\een
If we now put $k=N-M$ and apply this operator to the vector
\ben
e_{a_1}\ot
\cdots\ot\wh e_{a_i}\ot\cdots
\ot e_{a_{N-M}}\ot e_{a_j}
\een
then the coefficient at the vector
$
e_{a_1}\ot
\cdots\ot e_{a_{N-M}}
$
will be an expression involving only the entries $s_{ab}(u)$
of the matrix $S(u)$ with $a,b\in\Ac$. This proves the claim
for the auxiliary minors. The argument
for the Sklyanin minors is the same.

Finally, the calculation of the image of $s_{a_ia_j}(u)$
under $\vartheta$ is completed by the application of
$\varpi_{N}$ with the
use of \eqref{blockstw} and \eqref{blockstwsv} for the algebra
$\X(\g_N)$.
\epf

For any elements $a,b\in\Ac$ introduce the Sklyanin minors
\beql{wts}
s^{\sharp}_{ab}(u)={s}^{-m\cdots\ts m,\ts a}_{-m\cdots\ts m,\ts b}\ts
(u+M/2)
\eeq
and denote by $S^{\sharp}(u)$ the $(N-M)\times(N-M)$ matrix whose
$ab$-entry is $s^{\sharp}_{ab}(u)$.
Also, denote by $S^{}_{\Bc\Bc}(u)$ the submatrix of $S(u)$ whose rows
and columns are numbered by the elements of $\Bc$.

The following is a version of the {\it quantum Sylvester theorem\/}
for the extended twisted Yangian $\X(\g_N)$; cf. \cite{gr:dm},
\cite{kl:mi}, \cite{m:yt}.

\bth\label{thm:sylveven}
The mapping
\beql{homomsy}
s_{ab}(u)\mapsto s^{\sharp}_{ab}(u)
\eeq
defines an algebra homomorphism $\X(\g_{N-M})\to \X(\g_N)$. Moreover,
\begin{multline}
\sdet S^{\sharp}(u)=\sdet S(u+M/2)\\
{}\times{}\sdet S^{}_{\Bc\Bc}(u+M/2-1)
\cdots \ts\sdet S^{}_{\Bc\Bc}(u+M/2-N+M+1).
\non
\end{multline}
\eth

\bpf
By Proposition~\ref{prop:psimskmin}, we have
\beql{psimqmqstw}
\vartheta: s_{ab}(u)\mapsto
\big[\sdet S^{}_{\Bc\Bc}(u+M/2)\big]^{-1}\cdot
s^{\sharp}_{ab}(u).
\eeq
As we observed in the proof of Proposition~\ref{prop:psimskmin},
the expansion of
$\sdet S^{}_{\Bc\Bc}(u)$ only contains the series $s_{ij}(u)$
with $i,j\in\{-m,\dots,m\}$ (and with some shifts in $u$).
Hence, by Corollary~\ref{cor:centsm}, the series $\sdet S^{}_{\Bc\Bc}(u)$
commutes with $s^{\sharp}_{ab}(v)$ for any
$a,b\in\Ac$. Since $\vartheta$
preserves the defining relations
of the extended twisted Yangian,
we can conclude that the assignment \eqref{homomsy}
defines a homomorphism.

Furthermore, applying
\eqref{sdetcirc} and \eqref{blockstw},
we deduce
\ben
\vartheta:{s\ts}^{a_1\cdots\ts a_{N-M}}_{a_1\cdots\ts a_{N-M}}(u)
\mapsto
\big[\sdet S^{}_{\Bc\Bc}(u+M/2)\big]^{-1}\cdot\sdet S(u+M/2).
\een
Now we calculate this image in a different way and compare
the results. The expansion of the Sklyanin minor
${s\ts}^{a_1\cdots\ts a_{N-M}}_{a_1\cdots\ts a_{N-M}}(u)$
has the form of a linear combination
of the products
\ben
s_{\tss b_1c_1}(u)\ts s_{\tss b_2c_2}(u-1)\cdots
s_{\tss b_{N}c_{N}}(u-N+M+1),\qquad b_i,c_i\in\Ac,
\een
the coefficients being rational functions in $u$.
Using the relation \eqref{psimqmqstw}
we derive the desired
formula.
\epf

Denote by $\wh\si_{ab}(u)$ the entries of the Sklyanin comatrix
corresponding to the matrix $S^{\sharp}(u)$.

\bpr\label{prop:nnentry}
For any $a,b\in\Ac$ we have the relation
\begin{multline}
\wh\si_{ab}(u)=\wh s_{ab}(u+M/2)\\
{}\times{}\sdet S^{}_{\Bc\Bc}(u+M/2-1)
\cdots \ts\sdet S^{}_{\Bc\Bc}(u+M/2-N+M+2).
\non
\end{multline}
\epr

\bpf
We use the same argument as in the proof of Theorem~\ref{thm:sylveven}.
Let us calculate the image of the auxiliary minor
${\check s\ts}^{a_1\cdots\ts a_{N-M}}_{a_1\cdots\tss
\wh a_i\tss\cdots \ts a_{N-M},\ts a_j}(u)$
under the map
$\vartheta$ in two different ways; see \eqref{vartheta}.
Using \eqref{blockstwsv} we get
\begin{multline}
\varpi_{N-M}\big({\check s\ts}^{a_1\cdots\ts a_{N-M}}_{a_1\cdots\tss
\wh a_i\tss\cdots \ts a_{N-M},\ts a_j}(u)\big)\\
=(-1)^{N-M-i}\ts
\big[{s\ts}^{a_1\cdots\ts a_{N-M}}_{a_1\cdots\ts a_{N-M}}(-u+n-m-1)\big]^{-1}
\cdot s_{a_ia_j}(-u+n-m-1).
\non
\end{multline}
Now apply $\upsilon$ and $\varpi_{N}$
by using \eqref{blockstw} and \eqref{blockstwsv}
to conclude that
\ben
\vartheta:{\check s\ts}^{a_1\cdots\ts a_{N-M}}_{a_1\cdots\tss
\wh a_i\tss\cdots \ts a_{N-M},\ts a_j}(u)\mapsto
\big[\sdet S^{}_{\Bc\Bc}(u+M/2)\big]^{-1}\cdot
\wh s_{a_ia_j}(u+M/2).
\een
On the other hand,
the expansion of the auxiliary minor
${\check s\ts}^{a_1\cdots\ts a_{N-M}}_{a_1\cdots\tss
\wh a_i\tss\cdots \ts a_{N-M},\ts a_j}(u)$ in
terms of the matrix elements
has the form of a linear combination
of products
\ben
s_{b_1c_1}(u)\ts s_{b_2c_2}(u-1)\cdots
s_{b_{N-1}c_{N-1}}(u-N+2),\qquad b_i,c_i\in\{1,\dots,N\},
\een
the coefficients being rational functions in $u$.
The proof is completed by the application of
\eqref{psimqmqstw}.
\epf

Interchanging the roles of the sets $\Ac$ and $\Bc$ in the above arguments
one can easily derive the corresponding
dual versions of Theorem~\ref{thm:sylveven}
and Proposition~\ref{prop:nnentry}. Here we only record the counterpart
of the first part of Theorem~\ref{thm:sylveven} which will be used below.

\bpr\label{prop:dualhom}
The mapping
\ben
s_{ij}(u)\mapsto
{s}^{-n\ts\cdots\ts -m-1,\ts i,\ts m+1\cdots \ts n}_{-n\ts\cdots\ts -m-1,
\ts j,\ts m+1\cdots \ts n}\ts(u+n-m),\qquad -m\leqslant i,j\leqslant m
\een
defines an algebra homomorphism $\X(\g_{M})\to \X(\g_N)$.
\epr

\section{Sylvester theorem for the twisted Yangian}
\label{sec:ty}
\setcounter{equation}{0}

The {\it twisted Yangian\/} $\Y(\g_N)$ corresponding to the Lie algebra $\g_N$
is the quotient of the extended twisted Yangian $\X(\g_N)$ by
the following {\it symmetry relation\/}
\beql{symme}
\theta_{ij}s_{-j,-i}(-u)=s_{ij}(u)\pm
\frac{s_{ij}(u)-s_{ij}(-u)}{2u},
\eeq
or, in the matrix form,
\ben
S^t(-u)=S(u)\pm
\frac{S(u)-S(-u)}{2u}.
\een
From now on, we shall mainly work with the twisted Yangian and so
we keep the same notation $s_{ij}^{(r)}$ for the generators of the
algebra $\Y(\g_N)$.
Note that for any even series $g(u)\in 1+\CC[[u^{-2}]]\ts u^{-2}$
the mapping
\beql{autog}
s_{ij}(u)\mapsto g(u)\ts s_{ij}(u)
\eeq
defines an automorphism of $\Y(\g_N)$.

As we shall see below, the homomorphism
of Theorem~\ref{thm:sylveven} respects the symmetry relation
in the orthogonal case, while in the symplectic case a minor correction
is needed to obtain a corresponding homomorphism of the twisted Yangians.
In order to treat both cases simultaneously, introduce the following notation
\ben
\al_{p}(u)=\begin{cases}1\qquad&\text{in the orthogonal case}\\[4pt]
             \dfrac{u+1/2}{u-p+1/2}\qquad&\text{in the symplectic case}.
               \end{cases}
\een

The image of the Sklyanin determinant in the twisted Yangian $\Y(\g_N)$
acquires the following symmetry property
\beql{symsdet}
\al_{n}(u)^{-1}\cdot\sdet S(u)=\al_{n}(-u+N-1)^{-1}\cdot\sdet S(-u+N-1),
\eeq
see \cite[Section~4.11]{mno:yc}. Moreover, the mapping
\beql{comautotw}
S(u)\mapsto \al_{n}(u)\cdot\wh S(-u+N/2-1)
\eeq
defines a homomorphism of $\Y(\g_N)$ into itself;
see \cite[Proposition~2.1]{m:fd}.

We can now prove a quantum Sylvester theorem for the twisted Yangian $\Y(\g_N)$.
We use the notation of the previous section. In particular, recall that
$s^{\sharp}_{ab}(u)$ denotes the quantum minor as in \eqref{wts}
for any $a,b\in\Ac$.

\bth\label{thm:sylvtw}
The mapping
\beql{homomtw}
s_{ab}(u)\mapsto \al_{-m}(u)\ts s^{\sharp}_{ab}(u)
\eeq
defines an algebra homomorphism $\Y(\g_{N-M})\to \Y(\g_N)$.
Moreover,
\begin{multline}
\sdet \big[\al_{-m}(u)\ts S^{\sharp}(u)\big]=\al(u)\cdot\sdet S(u+M/2)\\
{}\times{}\sdet S^{}_{\Bc\Bc}(u+M/2-1)
\cdots \ts\sdet S^{}_{\Bc\Bc}(u+M/2-N+M+1),
\non
\end{multline}
where
\ben
\al(u)=\al_{-m}(u)\ts\al_{-m}(u-1)\cdots\al_{-m}(u-N+M+1).
\een
\eth

\bpf
Denote by $S^*(u)$ the matrix which occurs on the
right hand side of \eqref{comautotw}.
Then by \eqref{symsdet} we have
\ben
\varpi_{N}(S(u))=\frac{\al_n(u+N/2)}{c\ts(u+N/2)}\cdot S^*(u),
\een
where we have put $c(u)=\sdet S(u)$ for brevity.
Denote by $S^{*}_{\Ac\Ac}(u)$ the submatrix of $S^{*}(u)$
whose rows and columns are numbered by the elements of $\Ac$
and let $\wh s^{\ts*}_{ab}(u)$ denote the entries
of the Sklyanin comatrix corresponding to
$S^{*}_{\Ac\Ac}(u)$.
Using the definition of the auxiliary minors
and applying \eqref{comam} and \eqref{blockstwsv},
we come to the relation
\begin{multline}\label{sminstar}
{s}^{-m\ts\cdots\ts m,\ts a_i}_{-m\ts\cdots\ts m,\ts a_j}(u+M/2)=
c\ts(u+M/2)\\
{}\times
\prod_{i=1}^{N-M-1}\frac{\al_n(-u+n-m+N/2-i)}{c\ts(-u+n-m+N/2-i)}
\cdot\wh s^{\ts*}_{a_ia_j}(-u+n-m-1),
\end{multline}
where we have also used the fact that the coefficients of the
Sklyanin determinant are central in the twisted Yangian $\Y(\g_N)$.
Observe that by \eqref{symsdet} we have
\ben
c\ts(u+M/2)\cdot\frac{\al_n(-u+n-m+N/2-1)}{c\ts(-u+n-m+N/2-1)}=\al_n(u+M/2).
\een
Therefore, \eqref{sminstar} takes the form
\ben
{s}^{-m\ts\cdots\ts m,\ts a_i}_{-m\ts\cdots\ts m,\ts a_j}(u+M/2)=
\al_n(u+M/2)\cdot \varphi(u)
\cdot\wh s^{\ts*}_{a_ia_j}(-u+n-m-1),
\een
where
\ben
\varphi(u)=\prod_{i=2}^{N-M-1}\frac{\al_n(-u+n-m+N/2-i)}{c\ts(-u+n-m+N/2-i)}.
\een
By the symmetry property \eqref{symsdet} we have $\varphi(u)=\varphi(-u)$, and so
the multiplication of the
generator series $s_{ij}(u)$ by $\varphi(u)$ preserves
the twisted Yangian defining relations.
Furthermore,
by \eqref{comautotw}, the mapping
\ben
s_{ab}(u)\mapsto\al_{n-m}(u)\ts\wh s^{\ts*}_{ab}(-u+n-m-1),\qquad a,b\in\Ac
\een
preserves the defining relations of
 $\Y(\g_{N-M})$. Thus, we may conclude
that the mapping
\ben
s_{ab}(u)\mapsto \al_{n-m}(u)\ts\al_n(u+M/2)^{-1}\ts
{s}^{-m\ts\cdots\ts m,\ts a_i}_{-m\ts\cdots\ts m,\ts a_j}(u+M/2)
\een
defines a homomorphism $\Y(\g_{N-M})\mapsto \Y(\g_{N})$.
To complete the proof, observe that $\al_{n-m}(u)\ts\al_n(u+M/2)^{-1}=\al_{-m}(u)$.

The formula for the Sklyanin determinant of the matrix $\al_{-m}(u) S^{\sharp}(u)$
is immediate from Theorem~\ref{thm:sylveven} and
the definition of $\sdet S^{\sharp}(u)$.
\epf

The corresponding version of Proposition~\ref{prop:dualhom}
for the twisted Yangian has the following form.

\bpr\label{prop:dualhomtw}
The mapping
\ben
s_{ij}(u)\mapsto\al_{m-n}(u)\cdot
{s}^{-n\ts\cdots\ts -m-1,\ts i,\ts m+1\cdots \ts n}_{-n\ts\cdots\ts -m-1,
\ts j,\ts m+1\cdots \ts n}\ts(u+n-m),\qquad -m\leqslant i,j\leqslant m
\een
defines an algebra homomorphism $\Y(\g_{M})\to \Y(\g_N)$.\qed
\epr

Now we shall demonstrate that in the case of the twisted Yangian
the entries of the Sklyanin
comatrix can be expressed in terms of Sklyanin minors;
cf. \eqref{comam}. We need the following lemma.

\ble\label{lem:asrtas}
For the twisted Yangian $\Y(\g_N)$ we have
\beql{ansras}
A_N\ts S_1(u)\ts R^{\tss t}_{12}(-2u+1)\cdots
R^{\tss t}_{1N}(-2u+N-1)=\frac{2u+1}{2u\pm 1}
\ts A_N\ts {S}^{\tss t}_1(-u).
\eeq
\ele

\bpf
We have the relation
$(N-1)!\ts A_N=A_N\ts A^{\prime}_{N-1}$, where $A^{\prime}_{N-1}$
denotes the anti-symmetrizer corresponding
to the subset of indices $\{2,\dots,N\}$.
By \eqref{araq},
\ben
A^{\prime}_{N-1}\ts
R^{\tss t}_{12}(-2u+1)\cdots
R^{\tss t}_{1N}(-2u+N-1)=
A^{\prime}_{N-1}\ts\Big(1+\frac{Q_{12}+\cdots+Q_{1N}}{2u-1}\Big).
\een
Therefore, the left hand side of \eqref{ansras} takes the form
\ben
A_N\ts S_1(u)\ts \Big(1+\frac{Q_{12}+\cdots+Q_{1N}}{2u-1}\Big).
\een
Apply this operator to a basis vector
\ben
v_{ij}=e_{a_j}\otimes e_{a_1}\ot
\cdots\otimes \wh e_{a_i}\otimes\cdots
\otimes e_{a_{N}},\qquad i,j\in\{1,\dots,N\},
\een
where $(a_1,\dots,a_N)=(-n,\dots,n)$. The coefficient
at $A_N\ts v_{ii}$ will be equal to $s_{a_ia_j}(u)$ if $a_j=-a_i$, and
equal to the expression
\ben
\frac{2u}{2u-1}s_{a_ia_j}(u)\mp\frac{1}{2u-1}
\theta_{a_ia_j}s_{-a_j,-a_i}(u)
\een
if $a_j\ne -a_i$. In both cases the coefficient coincides with
\ben
\frac{2u+1}{2u\pm1}\ts s^t_{a_ia_j}(-u)
\een
due to the symmetry relation \eqref{symme}.
\epf

As before, we set $(a_1,\dots,a_N)=(-n,\dots,n)$.

\bpr\label{prop:skcom}
For any $i,j\in\{1,\dots,N\}$ we have the relation
\ben
\wh s_{a_ia_j}^{\ts t}(u)=
(-1)^{i+j}\cdot\al_{N-1}(u)\cdot
{s\ts}^{a_1\ts\cdots\ts\wh a_j\cdots\ts a_N}_{a_1\ts\cdots
\ts\wh a_i\cdots\ts a_N}(-u+N-2).
\een
\epr

\bpf
The relations \eqref{comdef} and \eqref{asinvnew} with $M=N-2$
imply that
\beql{asinvmn}
A_N \ts\langle S_1,\ldots,S_{N-1}\rangle
=A_N\ts\wh S_N(u)\ts (R^{\tss t}_{N-1,N})^{-1}
\cdots (R^{\tss t}_{1,N})^{-1}.
\eeq
Using the notation \eqref{ucircnew} and the definition
of the matrix $S^*(u)$ we can write $\wh S(u)=\al_n(u)\ts S^*(u^{\ts\circ}_N)$
so that \eqref{asinvmn} becomes
\beql{asinvmc}
A_N \ts\langle S_1,\ldots,S_{N-1}\rangle
=\al_n(u)\ts A_N\ts S^*_N(u^{\ts\circ}_N)\ts R^{\ts\circ}_{N-1,N}\cdots
R^{\ts\circ}_{1N}.
\eeq
However, $S(u)\mapsto S^*(u)$
defines a homomorphism $\Y(\g_N)\to \Y(\g_N)$. Therefore, writing
$A_N=\sgn\sigma \cdot A_N\ts P_{\sigma}$, where $\sigma=(1N)(2,N-1)\cdots$,
and applying Lemma~\ref{lem:asrtas} we can simplify
the right hand side of \eqref{asinvmc} as
\ben
\al_n(u)\ts A_N\ts S^*_N(u^{\ts\circ}_N)\ts R^{\ts\circ}_{N-1,N}\cdots
R^{\ts\circ}_{1N}=\al_n(u)\ts \frac{2u^{\ts\circ}_N+1}{2u^{\ts\circ}_N\pm 1}
\ts A_N\ts \ts {S^*}^t_N(-u^{\ts\circ}_N).
\een
Hence, we come to the identity
\ben
A_N \ts\langle S_1,\ldots,S_{N-1}\rangle
=\al_{n-1}(u)\ts A_N\ts \ts {S^*}^t_N(-u^{\ts\circ}_N).
\een
Applying both sides to the basis vector
$e_{a_1}\ot
\cdots\otimes \wh e_{a_i}\otimes\cdots
\otimes e_{a_{N}}\otimes e_{a_j}$ we get
\ben
\al_{n-1}(u)\ts s^{*\ts t}_{a_ia_j}(u-N/2+1)=(-1)^{i+j}\ts
{s\ts}^{a_1\ts\cdots\ts\wh a_j\cdots\ts a_N}_{a_1\ts\cdots
\ts\wh a_i\cdots\ts a_N}(u).
\een
The argument is completed by using the definition of $S^*(u)$.
\epf

Suppose now that $m=n-1$. As before, we identify
the subalgebra of $\Y(\g_N)$ generated by the elements
$s_{ab}^{(r)}$ for $a,b\in \{-n,n\}$
with the twisted Yangian $\Y(\g_2)$.
The restriction of the homomorphism \eqref{comautotw}
to the subalgebra $\Y(\g_2)$ defines the homomorphism
\ben
\phi:\Y(\g_2)\to\Y(\g_N).
\een
The defining relations of the twisted Yangian imply that the map
\ben
s_{ab}(u)\mapsto \sgn a\cdot\sgn b\cdot s_{ab}(u),\qquad a,b\in\{-n,n\}
\een
defines an automorphism $\psi:\Y(\g_2)\to\Y(\g_2)$.

\bco\label{cor:homcoin}
The homomorphism $\phi$
coincides with the homomorphism \eqref{homomtw} in the symplectic case,
while in the orthogonal case the homomorphism \eqref{homomtw}
coincides with the composition $\phi\circ\psi$.
\eco

\bpf
This is immediate from Proposition~\ref{prop:skcom} and
the skew-symmetry of the Sklyanin minors.
\epf

The symmetry relation \eqref{symme} allows one to obtain the
following expansion of the auxiliary minors (see \cite[Proposition~4.4]{m:sd}):
if $-b_1\in\{a_1,\dots,a_{k-1},c\}$
and $c\notin\{-b_2,\dots,-b_{k-1}\}$ then
\beql{auzerosf}
{\check s\ts}^{a_1\cdots\ts a_{k-1},\ts c}_{b_1\cdots\ts b_{k-1},\ts c}(u)=
\frac{2u+1}{2u\pm1}\ts\sum_{i=1}^{k-1}(-1)^{i-1}\ts s^t_{a_ib_1}(-u)
\ts {s\ts}^{a_1\cdots\ts\wh a_i\ts\cdots \ts a_{k-1}}_{b_2\ts\cdots\ts b_{k-1}}(u-1).
\eeq
Using this relation together with \eqref{sscheck} and \eqref{auzero}
one can derive explicit formulas for the Sklyanin determinant
and some Sklyanin minors.
The formulas use a special map
\beql{mapsym}
\om_N:\Sym_N\to \Sym_{N},\qquad p\mapsto p'
\eeq
from the symmetric group $\Sym_N$ into itself which is
defined by the following inductive procedure.
Given a set of positive integers
$c_1<\cdots<c_N$ we
regard $\Sym_N$ as the group of their permutations.
If $N=2$ we define $\om_2$ as the map $\Sym_2\to \Sym_{2}$
whose image is the identity permutation.
For $N>2$ define a map from the set of ordered pairs $(c_k,c_l)$
with $k\ne l$ into itself by the rule
\beql{ordpair}
\begin{alignedat}{2}
(c_k,c_l)&\mapsto (c_l,c_k),&&\qquad k,l<N,\\
(c_k,c_N)&\mapsto (c_{N-1},c_k),&&\qquad k<N-1,\\
(c_N,c_k)&\mapsto (c_k,c_{N-1}),&&\qquad k<N-1,\\
(c_{N-1},c_N)&\mapsto (c_{N-1},c_{N-2}),\\
(c_{N},c_{N-1})&\mapsto (c_{N-1},c_{N-2}).
\end{alignedat}
\end{equation}
Let $p=(p^{}_1,\dots,p^{}_N)$ be a permutation of the indices
$c_1,\dots,c_N$. Its image under
the map $\om_N$
is the permutation $p_{}^{\ts\prime}=
(p^{\ts\prime}_1,\dots,p^{\ts\prime}_{N-1},c_N)$, where the pair
$(p^{\ts\prime}_1,p^{\ts\prime}_{N-1})$ is the image
of the ordered pair $(p^{}_1,p^{}_N)$ under the map \eqref{ordpair}.
Then the pair $(p^{\ts\prime}_2,p^{\ts\prime}_{N-2})$ is found as
the image of $(p^{}_2,p^{}_{N-1})$
under the map \eqref{ordpair} which is defined on the set
of ordered pairs of elements obtained from $(c_1,\dots,c_N)$
by deleting $p^{}_1$ and $p^{}_N$. The procedure is completed in the same
manner by determining consequently
the pairs $(p^{\ts\prime}_i,p^{\ts\prime}_{N-i})$.

Now suppose that $M$ is a positive integer and $M=2m$ or $M=2m+1$.
For the proof of the following formula for
the Sklyanin minor in the twisted Yangian
$\Y(\g_N)$ see \cite{m:sd}.

\bpr\label{prop:minsf}
Suppose that $a_1,\dots,a_M,b_M$ are arbitrary indices from
the set $\{-n,\dots,n\}$.
We have
\begin{multline}
{s\ts}^{-a_1\ts\cdots \ts -a_{M}}_{a_1\ts\cdots\ts a_{M-1},b_{M}}(u)\\
{}=\al_m(u)\sum_{p\in\Sym_M}
\sgn p\tss p'\cdot s^{\ts t}_{-a_{p(1)},a_{p'(1)}}(-u)
\cdots
s^{\ts t}_{-a_{p(m)},a_{p'(m)}}(-u+m-1)\\
{}\times{}s^{}_{-a_{p(m+1)},a_{p'(m+1)}}(u-m)\cdots
s^{}_{-a_{p(M)},b_{p'(M)}}(u-M+1),
\non
\end{multline}
where the $s^t_{ij}(u)$ denote the entries of the matrix $S^t(u)$.
\qed
\epr

\bre\label{rem:comat}
Although this determinant-like formula does not apply to
arbitrary Sklyanin minors
${s\ts}^{a_1\cdots\ts a_{k}}_{b_1\cdots\ts b_{k}}(u)$,
it provides explicit formulas for the Sklyanin determinant
$\sdet S(u)$ and the matrix elements $\wh s_{ij}(u)$
of the Sklyanin comatrix; see Proposition~\ref{prop:skcom}.
\ere

For the use in Section~\ref{sec:cc}
we shall prove the following simple property of the map \eqref{mapsym}.

\ble\label{lem:probi}
The map
$
\Sym_N\to\Sym_N
$
defined by
$
p\mapsto
p\ts (p')^{-1}
$is bijective.
\ele

\bpf
Suppose that $p$ and $q$ are two elements of $\Sym_N$ such
that $p\ts (p^{\ts\prime})^{-1}=q\ts (q^{\ts\prime})^{-1}$. It
suffices to show that $p=q$. By definition of the map
$\om_N$ we have $p^{\ts\prime}_N=q^{\ts\prime}_N=N$ which
implies that $p^{}_N=q^{}_N$. Then, due to the formulas
\eqref{ordpair}, we have $p^{\ts\prime}_1=q^{\ts\prime}_1$. Hence,
$p^{}_1=q^{}_1$. Now, since the pairs $(p^{}_1,p^{}_N)$ and
$(q^{}_1,q^{}_N)$ coincide, so do their images under the map
\eqref{ordpair}. In particular,
$p^{\ts\prime}_{N-1}=q^{\ts\prime}_{N-1}$. This implies that
$p^{}_{N-1}=q^{}_{N-1}$ and the proof is completed by repeating
this argument for the pairs $(p^{}_{i+1},p^{}_{N-i})$ and
$(q^{}_{i+1},q^{}_{N-i})$ with $i=1,2,\dots$.
\epf

\section{Skew representations}
\label{sec:sr}
\setcounter{equation}{0}

As before, we suppose that $N=2n$ or $N=2n+1$ so that
\beql{oos}
\g_N=\oa_{2n+1},\quad \spa_{2n},\quad\text{or}\quad\oa_{2n}.
\end{equation}
The finite-dimensional irreducible representations of $\g_N$
are in a one-to-one correspondence with $n$-tuples
$\lambda=(\lambda_1,\dots,\lambda_n)$
where the numbers $\lambda_i$ satisfy the conditions
\ben
\lambda_i-\lambda_{i+1}\in \ZZ_+ \qquad\text{for}\quad i=1,\dots,n-1,
\een
and
\ben
\begin{split}
-2\ts \lambda_1&\in\ZZ_+ \qquad\text{for}\quad \g_N=\oa_{2n+1},\\
-\lambda_1&\in\ZZ_+ \qquad\text{for}\quad \g_N=\spa_{2n},\\
-\lambda_1-\lambda_2&\in\ZZ_+ \qquad\text{for}\quad \g_N=\oa_{2n}.
\end{split}
\een
Such an $n$-tuple $\lambda$ is called the
{\it highest weight\/}\footnote{In a more common notation,
the highest weight is the $n$-tuple $(-\la_n,\dots,-\la_1)$.
In particular, in the symplectic case this $n$-tuple is a partition.}
of the corresponding representation which
we shall denote by $V(\lambda)$.
It contains a unique, up to a constant factor, nonzero vector $\xi$
(the {\it highest vector\/}) such that
\beq
\begin{aligned}
F_{ii}\ts\xi&=\lambda_i\ts\xi\qquad
&\text{for}&\quad i=1,\dots,n,\\
F_{ij}\ts\xi&=0\qquad
&\text{for}&\quad -n\leqslant i<j\leqslant n.
\end{aligned}
\non
\end{equation}
Let $M$ be a nonnegative integer such that $N-M$ is even and positive.
So, $M=2m$ or $M=2m+1$ for some $m<n$.
We shall identify the Lie algebra $\g_M$ with the subalgebra
of $\g_N$ spanned by the elements $F_{ij}$ with
the indices satisfying $-m\leqslant i,j\leqslant m$.
Denote by $V(\lambda)^+$ the subspace of $\g_M$-highest vectors
in $V(\lambda)$:
\begin{equation}
V(\lambda)^+=\{\eta\in V(\lambda)\ |\ F_{ij}\ts \eta=0,
\qquad -m\leqslant i<j\leqslant m\}.
\non
\end{equation}
Given a $\g_M$-weight
$\mu=(\mu_1,\dots,\mu_m)$ we denote by $V(\lambda)^+_{\mu}$
the corresponding weight subspace in $V(\lambda)^+$:
\begin{equation}
V(\lambda)^+_{\mu}=\{\eta\in V(\lambda)^+\ |\ F_{ii}\ts\eta=
\mu_i\ts\eta,\qquad i=1,\dots,m\}.
\non
\end{equation}
We have a natural
vector space isomorphism
$V(\lambda)^+_{\mu}\cong\Hom_{\g_M}(V(\mu),V(\lambda))$.

For any $i,j\in\{-n,\dots,n\}$ introduce the series in $u^{-1}$
with coefficients in the universal enveloping algebra $\U(\g_N)$ by
\ben
f_{ij}(u)=\delta_{ij}+F_{ij}\Big(u\pm\frac12\Big)^{-1}.
\een
The mapping
\beql{evaltw}
\pi : s_{ij}(u)\mapsto f_{ij}(u)
\eeq
defines a surjective homomorphism $\Y(\g_N)\to\U(\g_N)$ called
the {\it evaluation homomorphism\/}; see \cite{o:ty} and
\cite[Proposition~3.11]{mno:yc}. Let $F(u)$ denote the $N\times N$ matrix
whose $ij$-th entry is the series $f_{ij}(u)$.
We may introduce the Sklyanin minors
${f\ts}^{a_1\cdots\ts a_k}_{b_1\cdots\ts b_k}(u)$
of this matrix as
the images of the corresponding minors of the matrix $S(u)$ with
respect to the evaluation homomorphism,
\ben
\pi:{s\ts}^{a_1\cdots\ts a_k}_{b_1\cdots\ts b_k}(u)
\mapsto
{f\ts}^{a_1\cdots\ts a_k}_{b_1\cdots\ts b_k}(u).
\een
By Theorem~\ref{thm:sylvtw}, we have a homomorphism
$\Y(\g_{N-M})\to\U(\g_N)$ given by
\beql{homuea}
\rho:s_{ab}(u)\mapsto \al_{-m}(u)\ts{f\ts}^{-m\cdots\ts m,\ts a}_{-m
\cdots\ts m,\ts b}\ts(u+M/2).
\eeq
Due to Corollary~\ref{cor:centsm},
the image of this homomorphism is contained in the centralizer
$\U(\g_N)^{\g_M}$
of the subalgebra $\g_M$ in the universal enveloping algebra $\U(\g_N)$.
On the other hand, the vector space
$V(\lambda)^+_{\mu}$ is obviously a representation
of $\U(\g_N)^{\g_M}$.
Thus,  $V(\lambda)^+_{\mu}$ becomes equipped
with the $\Y(\g_{N-M})$-module structure defined via
the homomorphism $\rho$. We call this module
the {\it skew representation\/} of $\Y(\g_{N-M})$.
In the particular case $M=0$ (with even $N$)
the skew representation is just
the {\it evaluation module\/} $V(\la)$ over $\Y(\g_{N})$
defined via the evaluation homomorphism \eqref{evaltw}.

The universal enveloping algebra $\U(\g_{N-M})$ can be identified
with a subalgebra of the twisted Yangian $\Y(\g_{N-M})$ via the embedding
$
F_{ab}\mapsto s_{ab}^{(1)},
$
see \cite[Proposition~3.12]{mno:yc}.
The elements $F_{ab}$ are stable under the composition
of this embedding with the homomorphism $\rho$. Indeed, if $a\ne b$
then this is verified directly from the definition of the Sklyanin minors
with the use of \eqref{artrtsi}. If $a=b$ then we may assume without
loss of generality that $m=n-1$ and $a=\pm n$. So, the claim
follows from Corollary~\ref{cor:homcoin} and \eqref{symsdet}.
In other words, the restriction of the $\Y(\g_{N-M})$-module
$V(\lambda)^+_{\mu}$ to the subalgebra $\U(\g_{N-M})$ coincides
with its natural action defined by the
$\U(\g_N)$-action on $V(\la)$.

We shall now concentrate on the symplectic case $\g_N=\spa_{N}$
where $N=2n$. Our next goal is to prove that the skew representations of
$\Y(\spa_{N-M})$ are irreducible.
First we show the following.

\bpr\label{prop:image}
The centralizer $\U(\spa_{N})^{\spa_{M}}$
is generated by the image of the homomorphism $\rho$
and the center of $\U(\spa_{N})$.
\epr

\bpf
For a different homomorphism
$\rho':\Y(\spa_{N-M})\to\U(\spa_{N})^{\spa_{M}}$
given by
\beql{rhoprmbr}
s_{ab}(u)\mapsto \al_{n}(u)\ts\wh f_{ab}(-u+n-1),\qquad a,b\in\Ac,
\eeq
where $\wh f_{ab}(u)$ is the image of $\wh s_{ab}(u)$ under
the evaluation homomorphism \eqref{evaltw},
this statement was proved in \cite[Section~4]{mo:cc}.
Denote by $\U'$ the subalgebra of $\U(\spa_{N})^{\spa_{M}}$
generated by the center of $\U(\spa_{N})$ and
the coefficients of the series
${f\ts}^{-m\cdots\ts m,\ts a}_{-m\cdots\ts m,\ts b}\ts(u)$
with $a,b\in\Ac$. It is sufficient to prove
that all the coefficients
of the series $\wh f_{ab}(u)$ belong to $\U'$.
However, the center of $\U(\spa_{N})$ is generated by the
coefficients of the Sklyanin determinant $\sdet F(u)$, i.e., the
image of $\sdet S(u)$ under the evaluation homomorphism \eqref{evaltw};
see \cite[Theorem~5.2]{m:sd}. Then by Theorem~\ref{thm:sylvtw}
the coefficients of the series $\sdet F_{\Bc\Bc}(u)$
belong to $\U'$.
Due to Proposition~\ref{prop:nnentry},
the coefficients
of the series $\wh f_{ab}(u)$ also belong to $\U'$.
\epf

\bco\label{cor:irred}
The skew representation $V(\lambda)^+_{\mu}$ of the twisted
Yangian $\Y(\spa_{N-M})$ is irreducible.
\eco

\bpf
Since the representation $V(\lambda)^+_{\mu}$ of
$\U(\spa_{N})^{\spa_{M}}$
is irreducible \cite[Section~9.1]{d:ae},
the statement follows from Proposition~\ref{prop:image} as
the central elements of $\U(\spa_{N})$ act on $V(\lambda)^+_{\mu}$
by scalar operators.
\epf

\bre\label{rem:ort} In general, the corresponding statement for the orthogonal
twisted Yangian $\Y(\oa_{N-M})$ is false; see \cite{m:gtb}
for the particular case $N-M=2$. If $N$ is even, then the $\Y(\oa_{2})$-module
$V(\lambda)^+_{\mu}$ is still irreducible. If $N$ is odd,
then for general parameters $\la$ and $\mu$
the $\Y(\oa_{2})$-module $V(\lambda)^+_{\mu}$ is isomorphic to the direct sum
of two irreducibles. It looks plausible that, in general,
the $\Y(\oa_{N-M})$-module
$V(\lambda)^+_{\mu}$ is completely reducible.
It would be interesting to obtain
its irreducible decomposition.
\ere

Now recall the classification results
for representations of the twisted Yangian
$\Y(\spa_{N})$; see \cite{m:fd}.
If $V$ is a finite-dimensional irreducible representation
of $\Y(\spa_{N})$ then $V$ contains a unique, up to a scalar
factor, vector $\xi\ne 0$ such that
\ben
\begin{aligned}
s_{ii}(u)\ts\xi&=\mu_i(u)\ts\xi\qquad
&\text{for}&\quad i=1,\dots,n,\\
s_{ij}(u)\ts\xi&=0\qquad
&\text{for}&\quad -n\leqslant i<j\leqslant n,
\end{aligned}
\een
where each $\mu_i(u)$ is a formal series in $u^{-1}$
with coefficients in $\CC$. Moreover, there exist monic
polynomials $P_1(u),\dots,P_n(u)$ in $u$
with $P_1(u)=P_1(-u+1)$
such that
\beql{dpk}
\frac{\mu_{i-1}(u)}{\mu_{i}(u)}=\frac{P_i(u+1)}{P_i(u)},\qquad i=2,\dots,n
\eeq
and
\beql{dpfirst}
\frac{\mu_{1}(-u)}{\mu_{1}(u)}=\frac{P_1(u+1)}{P_1(u)}.
\eeq
The $n$-tuple $\mu(u)=(\mu_1(u),\dots,\mu_n(u))$ is called
the {\it highest weight\/} and
the $P_i(u)$ are called the {\it Drinfeld polynomials\/}
of the representation $V$. Furthermore, given an
$n$-tuple $(P_1(u),\dots,P_n(u))$ of monic polynomials
with $P_1(u)=P_1(-u+1)$ there exists a
finite-dimensional irreducible representation $V$ of $\Y(\spa_{N})$
having this $n$-tuple as the family of its Drinfeld
polynomials. The isomorphism class of such representation
$V$ is determined uniquely,
up to the twisting by an automorphism of $\Y(\spa_{N})$
of the form
\ben
s_{ij}(u)\mapsto g(u)\ts s_{ij}(u),
\een
where $g(u)$ is a series in $u^{-2}$ with constant term $1$.

Since the Sklyanin determinant $\sdet S(u)$ is central in $\Y(\spa_{N})$,
it acts on $V$ by scalar multiplication. The scalar can be
calculated with the use of Proposition~\ref{prop:minsf}; see also \cite{m:sd}.
We have
\beql{sdetav}
\sdet S(u)|^{}_{V}=\al_n(u)\ts \prod_{i=1}^n \mu_i(-u+i-1)\ts\mu_i(u-N+i).
\eeq

\bpr\label{prop:sklmu}
With the above notation, for any $k=1,\dots,n$
in the representation $V$ we have
\ben
{s\ts}^{-k+1\ts\cdots\ts k}_{-k+1\ts\cdots\ts k}(u)\ts\xi=
\mu_k(u-2k+2)\ts
{s\ts}^{-k+1\ts\cdots\ts k-1}_{-k+1\ts\cdots\ts k-1}(u)\ts\xi.
\een
\epr

\bpf
By \eqref{sscheck} we can write
\ben
{s\ts}^{-k+1\ts\cdots\ts k}_{-k+1\ts\cdots\ts k}(u)=\sum_{c=-n}^n
{\check s\ts}^{-k+1\ts\cdots\ts k}_{-k+1\ts\cdots\ts k-1,c}(u)\ts
\ts s_{c\tss k}(u-2k+2).
\een
Apply both sides to the highest vector $\xi$. We have
$s_{c\tss k}(u-2k+2)\ts\xi=0$ if $c<k$.
On the other hand, if $c>k$ then
${\check s\ts}^{-k+1\ts\cdots\ts k}_{-k+1\ts\cdots\ts k-1,c}(u)=0$
by \eqref{auzero}. Finally, if $c=k$ then \eqref{aucc} gives
\ben
{\check s\ts}^{-k+1\ts\cdots\ts k}_{-k+1\ts\cdots\ts k-1,c}(u)=
{s\ts}^{-k+1\ts\cdots\ts k-1}_{-k+1\ts\cdots\ts k-1}(u)
\een
completing the proof.
\epf

Our aim now is to identify the representation
$V(\lambda)^+_{\mu}$ of $\Y(\spa_{N-M})$
by calculating its highest weight and Drinfeld polynomials.
Note that for the evaluation module $V(\la)$ over $\Y(\spa_{N})$
these can be immediately found from
\eqref{evaltw}. In particular, the $i$-th component of the
highest weight is given by
\beql{laimui}
\frac{u+\la_i-1/2}{u-1/2},\qquad i=1,\dots,n.
\eeq

In the case $M>0$
we employ the basis in the $\spa_{N}$-module $V(\lambda)$
constructed in \cite{m:br}. This basis is parameterized by
the {\it patterns\/} $\La$ which are arrays of non-positive integers
of the form
\begin{align}
\quad&\qquad\lambda^{}_{n1}\qquad\lambda^{}_{n2}
\qquad\qquad\cdots\qquad\qquad\lambda^{}_{nn}\non\\
&\lambda'_{n1}\qquad \lambda'_{n2}
\qquad\qquad\cdots\qquad\qquad\lambda'_{nn}\non\\
\ &\qquad\lambda^{}_{n-1,1}\qquad\cdots
\qquad\lambda^{}_{n-1,n-1}\non\\
&\lambda'_{n-1,1}
\qquad\cdots\qquad\lambda'_{n-1,n-1}\non\\
&\qquad\cdots\qquad\cdots\non\\
\quad&\qquad\lambda^{}_{11}\non\\
&\lambda'_{11}\non
\end{align}
where $\lambda=(\lambda^{}_{n1},\dots, \lambda^{}_{nn})$
is the top row of $\La$
and the following {\it betweenness conditions\/} hold
\beq
0\geqslant\lambda'_{k1}\geqslant\lambda^{}_{k1}\geqslant\lambda'_{k2}\geqslant
\lambda^{}_{k2}\geqslant \cdots\geqslant
\lambda'_{k,k-1}\geqslant\lambda^{}_{k,k-1}\geqslant
\lambda'_{kk}\geqslant\lambda^{}_{kk}
\non
\end{equation}
for $k=1,\dots,n$, and
\beq
0\geqslant\lambda'_{k1}\geqslant\lambda^{}_{k-1,1}\geqslant\lambda'_{k2}\geqslant
\lambda^{}_{k-1,2}\geqslant \cdots\geqslant
\lambda'_{k,k-1}\geqslant\lambda^{}_{k-1,k-1}\geqslant\lambda'_{kk}
\non
\end{equation}
for $k=2,\dots,n$.
The representation $V(\la)$ admits a basis $\zeta^{}_{\Lambda}$
parameterized by all patterns $\La$. The formulas for the matrix
elements of a family of generators of the Lie algebra $\spa_N$
in this basis can be explicitly written down; see \cite{m:br}.
We shall only need these formulas for the generators $F_{kk}$.
We have
\beql{fkk}
F_{kk}\ts \zeta^{}_{\Lambda}=\left(2\sum_{i=1}^k\lambda'_{ki}-
\sum_{i=1}^k\lambda^{}_{ki}-\sum_{i=1}^{k-1}\lambda^{}_{k-1,i}\right)
\zeta^{}_{\Lambda},\qquad k=1,\dots,n.
\eeq

Given a highest weight $\mu=(\mu_1,\dots,\mu_m)$ for $\spa_M$,
a basis of the space $V(\lambda)^+_{\mu}$ is formed
by those vectors $\zeta^{}_{\La}$ for which the row
$(\la^{}_{m1},\dots,\la^{}_{mm})$ of $\La$ coincides with $\mu$
and
\ben
\la^{}_{ki}=\la'_{ki}=\mu^{}_i, \qquad 1\leqslant i\leqslant k\leqslant m.
\een
Omitting the same triangle part below the $m$-th row of all such
patterns $\La$
we get trapezium-like patterns
(still denoted by $\La$)
with the top row $\lambda$ and the bottom
row $\mu$, as illustrated:
\begin{align}
\quad&\qquad\lambda^{}_{1}\qquad\lambda^{}_{2}
\qquad\qquad\qquad\cdots\qquad\qquad\qquad\lambda^{}_{n}\non\\
&\lambda'_{n1}\qquad \lambda'_{n2}
\qquad\qquad\qquad\cdots\qquad\qquad\qquad\lambda'_{nn}\non\\
&\qquad\cdots\qquad\cdots\non\\
&\lambda'_{m+1,1}\quad \lambda'_{m+1,2}
\qquad\quad\cdots\quad\qquad\lambda'_{m+1,m+1}\non\\
\quad&\ \ \qquad\mu^{}_{1}\qquad\ \mu^{}_{2}
\qquad\cdots\qquad\mu^{}_{m}\non
\end{align}
Due to the betweenness conditions, the space $V(\lambda)^+_{\mu}$
is nonzero if and only if
\ben
\mu_i\geqslant\la_{i+n-m},\qquad i=1,\dots,m,
\een
and
\ben
\la_i\geqslant\mu_{i+n-m},\qquad i=1,\dots,n,
\een
assuming $\mu_i=-\infty$ for $i>m$. We also set
$\mu_i=0$ for $i\leqslant 0$.
In what follows we suppose that $V(\lambda)^+_{\mu}$
is nonzero.

Denote by $\h$ the diagonal Cartan subalgebra
of $\spa_{N-M}$ spanned by the basis vectors
$F_{kk}$ with $k=m+1,\dots,n$.
Let the $\ve_k\in\h^*$ with $k=m+1,\dots,n$ be the dual
basis vectors of $\h^*$.
Consider the root system for $\spa_{N-M}$ with respect to $\h$
where the basis vectors $F_{ij}$ with $i<j$
are positive root vectors so that the positive roots are
$-2\ve_{k}$ with $k=m+1,\dots,n$ and ${}\pm\ve_i-\ve_j$
with $m+1\leqslant i<j\leqslant n$.
By \eqref{fkk}, the weight $w(\La)=(w_{m+1},\dots,w_n)$
of a trapezium pattern $\La$ with respect to $\h$
is given by
\ben
w_k=2\sum_{i=1}^k\lambda'_{ki}-
\sum_{i=1}^k\lambda^{}_{ki}-\sum_{i=1}^{k-1}\lambda^{}_{k-1,i},
\qquad k=m+1,\dots,n.
\een
We shall need the standard partial ordering on the set of weights
of an $\spa_{N-M}$-module $V$. We shall write $w\preccurlyeq w'$
for two weights $w$ and $w'$ of $V$ if $w'-w$ is a linear
combination of positive roots with nonnegative integral coefficients.
Since the vectors $\zeta^{}_{\La}$ corresponding to the
trapezium patterns $\La$ form a basis of
the $\spa_{N-M}$-module $V(\lambda)^+_{\mu}$,
the set of weights of this module is comprised
by the weights $w(\La)$ for all possible patterns $\La$.

Given three integers $i,j,k$
we shall denote by $\middle\{i,j,k\}$ that of the three which is
between the two others. If one of the indices, say, $k$ is the symbol $-\infty$
then $\middle\{i,j,k\}$ is understood as $\min\{i,j\}$.
Consider the trapezium array $\La_0$ whose entries are determined
by
\beql{laki}
\la^{}_{ki}=\middle\{\la_i,\mu_{i+k-m},\mu_{i+m-k}\}
\eeq
and
\ben
\la'_{ki}=\middle\{\la_i,\mu_{i+k-m-1},\mu_{i+m-k}\}
\een
for all possible values of $i$ and $k$. One easily verifies
that $\La_0$ is a pattern.

\bpr\label{prop:umaxw}
The $\spa_{N-M}$-module $V(\lambda)^+_{\mu}$
has a unique maximal weight. This weight
coincides with $w(\La_0)$.
\epr

\bpf
Suppose that $w(\La)$ is a maximal weight
of $V(\lambda)^+_{\mu}$
for some pattern $\La$.
Then the entries of $\La$ should satisfy
\beql{lapr}
\la'_{ki}=\max\{\la^{}_{ki},\la^{}_{k-1,i}\}
\eeq
for all $k=m+1,\dots,n-1$ and $i=1,\dots,k$, where we assume
$\la^{}_{ki}=-\infty$ for $i>k$.
Indeed, if the equality is not attained for some $k$ and $i$
then by the betweenness conditions we have
$\la'_{ki}>\max\{\la^{}_{ki},\la^{}_{k-1,i}\}$.
Therefore,
decreasing the entry $\la'_{ki}$ by $1$ we get a pattern of greater
weight than $w(\La)$ which contradicts the maximality of $w(\La)$.
Thus, the weight $w(\La)=(w_{m+1},\dots,w_n)$ can now be
written as
\ben
w_k=\sum_{i=1}^{k-1}|\la^{}_{ki}-\la^{}_{k-1,i}|+\la^{}_{kk},
\qquad k=m+1,\dots,n.
\een
We shall argue by induction on $n-m$ to show that
the $\la^{}_{ki}$ must be given by \eqref{laki}. This will
imply the statement since the entries of $\La_0$
do satisfy \eqref{lapr} which is easily seen.
In the case $n-m=1$ there is nothing to prove. Suppose that $n-m=2$.
Omitting the primed entries, we can depict $\La$ as
\begin{align}
\quad&\qquad\lambda^{}_{1}\quad\lambda^{}_{2}
\quad\cdots\quad\la^{}_{n-2}\quad\la^{}_{n-1}\quad\la^{}_{n}\non\\
\quad&\qquad\rho^{}_{1}\quad\rho^{}_{2}
\quad\cdots\quad\rho^{}_{n-2}\quad\rho^{}_{n-1}\non\\
\quad&\qquad\mu^{}_{1}\quad\mu^{}_{2}
\quad\cdots\quad\mu^{}_{n-2}\non
\end{align}
where we have put $\rho_i=\la^{}_{n-1,i}$. The weight $w(\La)=(w_{n-1},w_n)$
is given by
\beql{weig}
\bal
w_{n-1}&=\sum_{i=1}^{n-2}|\rho_i-\mu_i|+\rho_{n-1},\\
w_{n}&=\sum_{i=1}^{n-1}|\la_i-\rho_i|+\la_{n}.
\eal
\eeq
Suppose that \eqref{laki} is violated for some $i$ so that
$\rho_i\ne\middle\{\la_i,\mu_{i+1},\mu_{i-1}\}$. Observe
that by the betweenness conditions we have
\ben
\rho_{i+1}\leqslant\middle\{\la_i,\mu_{i+1},\mu_{i-1}\}\leqslant
\rho_{i-1}.
\een
Hence, if $\rho_i<\middle\{\la_i,\mu_{i+1},\mu_{i-1}\}$
(respectively, $\rho_i>\middle\{\la_i,\mu_{i+1},\mu_{i-1}\}$) then
we can increase (respectively, decrease) the value of $\rho_i$ by $1$
without violating the betweenness conditions and thus
to get another pattern $\La'$. Due to the formulas \eqref{weig},
we have $w(\La)\prec w(\La')$ which
contradicts the maximality of $w(\La)$. This proves the statement
in the case under consideration.

Suppose now that $n-m>2$. Let us set $\rho_i=\la^{}_{n-1,i}$ as above,
and consider the set of patterns having $\rho=(\rho_1,\dots,\rho_{n-1})$
as the top row and $\mu$ as the bottom row.
Then, by the maximality of $w(\La)$, the trapezium subpattern of $\La$
with the top row $\rho$
and bottom row $\mu$ will clearly be of a maximal weight
amongst all patterns of this set.
By the induction hypothesis,
we must have
\beql{lakirho}
\la^{}_{ki}=\middle\{\rho_i,\mu_{i+k-m},\mu_{i+m-k}\}
\eeq
for all $k=m+1,\dots,n-2$ and $i=1,\dots,k$.
Similarly,
the subpattern of $\La$ having the top row $\la$ and the bottom row
$\si=(\si_1,\dots,\si_{n-2})$ is of a maximal weight amongst
all patterns of this form, where we have put $\si_i=\la^{}_{n-2,i}$.
By the statement for the case $n-m=2$, we must have
\ben
\rho_i=\middle\{\la_i,\si_{i+1},\si_{i-1}\},\qquad i=1,\dots,n-1.
\een
Combining this with the relations
$
\si_i=\middle\{\rho_i,\mu_{i+n-m-2},\mu_{i+m-n+2}\}
$
implied by \eqref{lakirho}, we get the desired relation
for the row $\rho$ of $\La$,
\ben
\rho_i=\middle\{\la_i,\mu_{i+n-m-1},\mu_{i+m-n+1}\},\qquad i=1,\dots,n-1.
\een
Indeed, this is easily verified by looking at all
possible combinations of the three values of each of $\si_{i-1}$ and $\si_{i+1}$.
Finally, substituting these values of $\rho_i$ into
\eqref{lakirho} we conclude that \eqref{laki}
holds for all possible $k$ and $i$. Thus, $\La$ coincides with $\La_0$.
\epf

\bco\label{cor:tyhv}
The vector $\zeta^{}_{\La_0}$ is the highest vector of the
$\Y(\spa_{N-M})$-module $V(\lambda)^+_{\mu}$.
\eco

\bpf
It suffices to demonstrate that $\zeta^{}_{\La_0}$ is annihilated
by all generators $s_{ab}(u)$ of $\Y(\spa_{N-M})$ with $a<b$,
since by \cite[Remark~4.4]{m:fd} the vector $\zeta^{}_{\La_0}$
will then have to be an eigenvector for all $s_{aa}(u)$ with
$a=m+1,\dots,n$.

Recall that $s_{ab}(u)$ acts on $V(\lambda)^+_{\mu}$ as
a Sklyanin minor which is a series in $u^{-1}$ with coefficients
in the universal enveloping algebra $\U(\spa_{N})$;
see \eqref{homuea}.
Using Proposition~\ref{prop:comsm},
it is easy to derive that the weight of all these coefficients coincides
with the weight of the element $F_{ab}$ with respect to the adjoint action
of the Cartan subalgebra $\h$ of $\spa_{N-M}$. Therefore, if $a<b$ then
the vector $s_{ab}(u)\ts\zeta^{}_{\La_0}$
has the weight $w(\La_0)+\al$ for a positive root $\al$.
Now the claim follows from Proposition~\ref{prop:umaxw}.
\epf

Corollary~\ref{cor:tyhv} implies that the highest weight
$\mu(u)=(\mu_{m+1}(u),\dots,\mu_n(u))$ of the
$\Y(\spa_{N-M})$-module $V(\lambda)^+_{\mu}$ is determined by
the relations
\ben
s_{aa}(u)\ts \zeta^{}_{\La_0}=\mu_{a}(u)\ts \zeta^{}_{\La_0},
\qquad a=m+1,\dots,n.
\een
In order to calculate $\mu(u)$ we shall use the results of \cite{m:br},
where the particular case $m=n-1$ was considered.
In \cite{m:br} the vector space $V(\lambda)^+_{\mu}$ was endowed with
the $\Y(\spa_{2})$-module structure defined by the composition
of an automorphism of the type \eqref{autog} and the homomorphism
\eqref{rhoprmbr}.
Therefore, using Corollary~\ref{cor:homcoin}, we can reformulate
the result for the $\Y(\spa_{2})$-module structure on $V(\lambda)^+_{\mu}$
defined by the homomorphism \eqref{homuea} to obtain the following.

\bpr\label{prop:mbr}
For $m=n-1$ in the
$\Y(\spa_{2})$-module $V(\lambda)^+_{\mu}$ we have
\begin{multline}
\al_{-n+1}(u)\ts {f\ts}^{-n+1\ts \cdots \ts n}_{-n+1\ts
\cdots \ts n}\ts(u+n-1)\ts \zeta^{}_{\La_0}\\
{}=\prod_{i=2}^n\frac{u-\min\{\la_{i-1},\mu_{i-1}\}+i-1/2}{u+i-1/2}
\cdot\prod_{i=1}^n\frac{u+\max\{\la_{i},\mu_{i}\}-i+1/2}{u-i+1/2}
\ts \zeta^{}_{\La_0}.
\non
\end{multline}
\epr

In the case of arbitrary $m<n$ introduce the following notation:
\beql{nuu}
\nu(u)=\prod_{i=1}^{m}\frac{(u+\mu_i-i+1/2)\ts(u-\mu_i+i+1/2)}
{(u-i+1/2)\ts(u+i+1/2)}.
\eeq

\bth\label{thm:hwmu}
The highest weight $\mu(u)=(\mu_{m+1}(u),\dots,\mu_n(u))$
of the $\Y(\spa_{N-M})$-module $V(\lambda)^+_{\mu}$
is given by the formulas
\ben
\bal
\mu_k(u)=\nu(u)\cdot\prod_{\overset{\scl{i=1}}{\la_i<\mu_{i+k-m-1}}}^{k-1}
&\frac{u-\max\{\la_i,\mu_{i+k-m}\}+k-m+i-1/2}{u-\mu_{i+k-m-1}+k-m+i-1/2}\\
{}\times\prod_{\overset{\scl{i=1}}{\la_i>\mu_{i+m-k+1}}}^{k-1}
&\frac{u+\min\{\la_i,\mu_{i+m-k}\}+k-m-i-1/2}{u+\mu_{i+m-k+1}+k-m-i-1/2}\\
{}\times{}
&\frac{u+\min\{\la_k,\mu_{m}\}-m-1/2}{u-m-1/2},
\eal
\een
where $k=m+1,\dots,n$.
\eth

\bpf
Let us denote by $\vp_{ab}(u)$ the image of the series $s_{ab}(u)$
under the homomorphism $\rho$ defined in \eqref{homuea}, that is,
\beql{phiab}
\vp_{ab}(u)= \al_{-m}(u)\ts{f\ts}^{-m\cdots\ts m,\ts a}_{-m
\cdots\ts m,\ts b}\ts(u+M/2).
\eeq
Combine these series into the $(N-M)\times(N-M)$ matrix $\Phi(u)$.
We shall use the usual notation for the
Sklyanin comatrix and the Sklyanin minors of $\Phi(u)$.
Applying Proposition~\ref{prop:nnentry}, we obtain the following
expression for the $nn$-th entry of the Sklyanin comatrix $\wh\Phi(u)$,
\beql{phif}
\wh\vp_{nn}(u)=\al_{-m}(u)\ts
\wh f_{nn}(u+M/2)\cdot f_{\Bc}(u-1)
\cdots \ts f_{\Bc}(u-N+M+2),
\eeq
where we have put
$f_{\Bc}(u)=\al_{-m}(u)\ts\sdet F^{}_{\Bc\Bc}(u+M/2)$.
On the other hand, by Proposition~\ref{prop:skcom}, we have
\ben
\wh\vp_{nn}(u)=\al_{N-M-1}(u)\ts
{\vp\ts}^{-n+1\ts \cdots\ts -m-1,\ts m+1\cdots \ts n}_{-n+1\ts
\cdots\ts -m-1,\ts m+1
\cdots \ts n}\ts(-u+N-M-2)
\een
and
\ben
\wh f_{nn}(u)=\al_{N-1}(u)\ts
{f\ts}^{-n+1\ts \cdots \ts n}_{-n+1\ts \cdots \ts n}\ts(-u+N-2).
\een
Therefore, by \eqref{phif},
\begin{multline}\label{phifsm}
\al_{N-M-1}(u)\ts
{\vp\ts}^{-n+1\ts \cdots\ts -m-1,\ts m+1\cdots \ts n}_{-n+1\ts \cdots\ts -m-1,\ts m+1
\cdots \ts n}\ts(-u+N-M-2)=\\
\al_{N-1}(u+M/2)\ts
{f\ts}^{-n+1\ts \cdots \ts n}_{-n+1\ts \cdots \ts n}\ts(-u+N-M/2-2)\\
{}\times{}\al_{-m}(u)\ts f_{\Bc}(u-1)
\cdots \ts f_{\Bc}(u-N+M+2).
\end{multline}
Similarly, applying Theorem~\ref{thm:sylvtw}, we derive the identity
\begin{multline}
\al_{N-M}(u)\ts
{\vp\ts}^{-n\ts \cdots\ts -m-1,\ts m+1\cdots \ts n}_{-n\ts \cdots\ts -m-1,\ts m+1
\cdots \ts n}\ts(-u+N-M-1)=\\
\al_{N}(u+M/2)\ts
{f\ts}^{-n\ts \cdots \ts n}_{-n\ts \cdots \ts n}\ts(-u+N-M/2-1)\\
{}\times{}\al_{-m}(u)\ts f_{\Bc}(u-1)
\cdots \ts f_{\Bc}(u-N+M+1).
\non
\end{multline}
Replacing here $n$ by $n-1$ and $u$ by $u-1$ we get
\begin{multline}\label{phifsdet}
\al_{N-M-2}(u-1)\ts
{\vp\ts}^{-n+1\ts \cdots\ts -m-1,\ts m+1\cdots \ts n-1}_{-n+1\ts
\cdots\ts -m-1,\ts m+1\cdots \ts n-1}\ts(-u+N-M-2)=\\
\al_{N-2}(u+M/2-1)\ts
{f\ts}^{-n+1\ts \cdots \ts n-1}_{-n+1\ts \cdots \ts n-1}\ts(-u+N-M/2-2)\\
{}\times{}\al_{-m}(u-1)\ts f_{\Bc}(u-2)
\cdots \ts f_{\Bc}(u-N+M+2).
\end{multline}
Due to Proposition~\ref{prop:sklmu}, in
$V(\lambda)^+_{\mu}$ we have
\begin{multline}
{\vp\ts}^{-n+1\ts \cdots\ts -m-1,\ts m+1\cdots \ts n}_{-n+1\ts \cdots\ts -m-1,\ts m+1
\cdots \ts n}\ts(v)
\ts\zeta^{}_{\La_0}\\
{}=\mu_n(v-N+M+2)\ts
{\vp\ts}^{-n+1\ts \cdots\ts -m-1,\ts m+1\cdots \ts n-1}_{-n+1\ts
\cdots\ts -m-1,\ts m+1\cdots \ts n-1}\ts(v)\ts\zeta^{}_{\La_0}.
\non
\end{multline}
Hence, comparing \eqref{phifsm} and \eqref{phifsdet}
we come to the relation
\begin{multline}
f_{\Bc}(u-1)\ts
{f\ts}^{-n+1\ts \cdots \ts n}_{-n+1\ts \cdots \ts n}\ts(-u+2n-m-2)
\ts\zeta^{}_{\La_0}\\
{}=\mu_n(-u)\ts
{f\ts}^{-n+1\ts \cdots \ts n-1}_{-n+1\ts \cdots \ts n-1}\ts(-u+2n-m-2)
\ts\zeta^{}_{\La_0}.
\non
\end{multline}
Therefore, for each $k=m+1,\dots,n$ the corresponding component
of the highest weight of the
$\Y(\spa_{N-M})$-module $V(\lambda)^+_{\mu}$ can be found
from the relation
\begin{multline}\label{muk}
f_{\Bc}(-u-1)\ts
{f\ts}^{-k+1\ts \cdots \ts k}_{-k+1\ts \cdots \ts k}\ts(u+2k-m-2)
\ts\zeta^{}_{\La_0}\\
{}=\mu_k(u)\ts
{f\ts}^{-k+1\ts \cdots \ts k-1}_{-k+1\ts \cdots \ts k-1}\ts(u+2k-m-2)
\ts\zeta^{}_{\La_0}.
\end{multline}
Observe that each of the Sklyanin minors
${f\ts}^{-k+1\ts \cdots \ts k}_{-k+1\ts \cdots \ts k}\ts(v)$
and
${f\ts}^{-k+1\ts \cdots \ts k-1}_{-k+1\ts \cdots \ts k-1}\ts(v)$
commutes with all elements of the subalgebra $\spa_{2k-2}$
by Corollary~\ref{cor:centsm}. Therefore,
since the basis $\{\zeta^{}_{\La}\}$ of $V(\la)$ is consistent
with the embeddings $\spa_{2k-2}\subset\spa_{2k}$ \cite{m:br}, the vector
$\zeta^{}_{\La_0}$
is an eigenvector for each of these minors. The corresponding
eigenvalue for the first minor can be calculated from
Proposition~\ref{prop:mbr}, so that
\beql{alfkk}
\bal
\al_{-k+1}(u)\ts &{f\ts}^{-k+1\ts \cdots \ts k}_{-k+1\ts
\cdots \ts k}\ts(u+k-1)\ts \zeta^{}_{\La_0}\\
{}&{}=\prod_{i=2}^k\frac{u-\min\{\la^{}_{k,i-1},\la^{}_{k-1,i-1}\}+i-1/2}{u+i-1/2}\\
{}&{}\times\prod_{i=1}^k\frac{u+\max\{\la^{}_{ki},\la^{}_{k-1,i}\}-i+1/2}{u-i+1/2}
\ts \zeta^{}_{\La_0},
\eal
\eeq
where the $\la_{ki}$ are the entries of $\La_0$. The second minor
coincides with the image of the Sklyanin determinant
for $\Y(\spa_{2k-2})$ under the evaluation homomorphism \eqref{evaltw}.
Hence the corresponding eigenvalue can be found from \eqref{sdetav} and
\eqref{laimui} which gives
\begin{multline}
\al_{-k+1}(u)\ts {f\ts}^{-k+1\ts \cdots \ts k-1}_{-k+1\ts
\cdots \ts k-1}\ts(u+k-1)\ts \zeta^{}_{\La_0}\\
{}=\prod_{i=2}^{k}\frac{u-\la^{}_{k-1,i-1}+i-1/2}{u+i-1/2}
\cdot\prod_{i=1}^{k-1}\frac{u+\la^{}_{k-1,i}-i+1/2}{u-i+1/2}
\ts \zeta^{}_{\La_0};
\non
\end{multline}
see also \cite{m:sd}.
Note that since
$f_{\Bc}(u)=\al_{-m}(u)\ts {f\ts}^{-m\ts \cdots \ts m}_{-m\ts
\cdots \ts m}\ts(u+m)$,
the last formula with $k=m+1$ also applies for the calculation
of $f_{\Bc}(u)\ts \zeta^{}_{\La_0}$. Hence,
\beql{fbhc}
f_{\Bc}(-u-1)\ts \zeta^{}_{\La_0}=\nu(u)\ts \zeta^{}_{\La_0}
\eeq
with $\nu(u)$ defined in \eqref{nuu}.
Furthermore, the formulas \eqref{laki}
for the entries of $\La_0$ imply that
\beql{maxmid}
\max\{\la^{}_{ki},\la^{}_{k-1,i}\}=\middle\{\la_i,\mu_{i+k-m-1},\mu_{i+m-k}\}
\eeq
and
\beql{minmid}
\min\{\la^{}_{ki},\la^{}_{k-1,i}\}=\middle\{\la_i,\mu_{i+k-m},\mu_{i+m-k+1}\}.
\eeq
Using \eqref{muk} we obtain the following expression for $\mu_k(u)$:
\begin{multline}
\mu_k(u)=\nu(u)\ts
\prod_{i=2}^{k}\frac{u-\middle\{\la_{i-1},\mu_{i+k-m-1},\mu_{i+m-k}\}+k-m+i-3/2}
{u-\middle\{\la_{i-1},\mu_{i+k-m-2},\mu_{i+m-k}\}+k-m+i-3/2}\\
{}\times
\prod_{i=1}^{k-1}\frac{u+\middle\{\la_{i},\mu_{i+k-m-1},\mu_{i+m-k}\}+k-m-i-1/2}
{u+\middle\{\la_{i},\mu_{i+k-m-1},\mu_{i+m-k+1}\}+k-m-i-1/2}\\
{}\times\frac{u+\min\{\la_k,\mu_m\}-m-1/2}{u-m-1/2}.
\non
\end{multline}
Finally, replace the index $i$ in the first product by $i+1$ and note
that if $\la_i\geqslant \mu_{i+k-m-1}$ then the corresponding factor
equals $1$. Similarly, if $\la_i\leqslant \mu_{i+m-k+1}$
then the corresponding factor in the second product
equals $1$. This brings the expression for $\mu_k(u)$
to the required form.
\epf

We can now
compute the Drinfeld polynomials for the
$\Y(\spa_{N-M})$-module
$V(\lambda)^+_{\mu}$
by using Theorem~\ref{thm:hwmu}.
Given any $\spa_{2n}$-highest weight $\la=(\la_1,\dots,\la_n)$,
set $\la_{-i}=-\la_i$ for $i=1,\dots, n$. We also assume that
$\la_0=0$ while $\la_k=-\infty$ and $\la_{-k}=+\infty$ for $k>n$.
Introduce the {\it diagram\/} $\Ga(\la)$ as a certain infinite set of unit squares
(cells) on the plane whose centers have integer coordinates.
The coordinates $(i,j)$ of a cell are interpreted as the row and column
number so that $i$ increases from top to bottom and $j$ increases
from left to the right.
With these assumptions\footnote{This definition of $\Ga(\la)$
corresponds to the one outlined in the Introduction
for the partition $(\la_{-n},\dots,\la_{-1})$.},
\ben
\Ga(\la)=\{(i,j)\in\ZZ^2\ |\ {-n}\leqslant i\leqslant n+1,\quad
\la_i\leqslant j< \la_{i-1}\}.
\een
The diagram has a central symmetry,
as illustrated below for $\la=(-4,-7)$ and $n=2$:

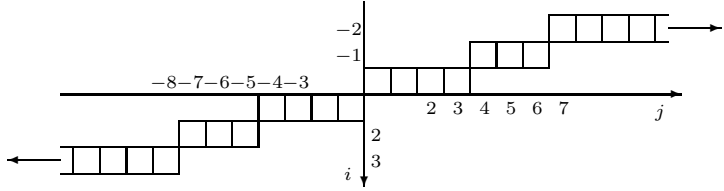
\setlength{\unitlength}{0.5em}
\begin{center}
\begin{picture}(48,14)

\put(1,0){\line(1,0){9}}
\put(1,2){\line(1,0){15}}
\put(10,4){\line(1,0){14}}
\put(16,4){\line(1,0){8}}
\put(16,6){\line(1,0){16}}
\put(24,6){\line(1,0){8}}
\put(24,8){\line(1,0){14}}
\put(32,8){\line(1,0){6}}
\put(32,10){\line(1,0){15}}
\put(38,12){\line(1,0){9}}

\multiput(2,0)(2,0){5}{\line(0,1){2}}
\multiput(10,2)(2,0){4}{\line(0,1){2}}
\multiput(16,4)(2,0){5}{\line(0,1){2}}
\multiput(24,6)(2,0){5}{\line(0,1){2}}
\multiput(32,8)(2,0){4}{\line(0,1){2}}
\multiput(38,10)(2,0){5}{\line(0,1){2}}

\put(1,1){\vector(-1,0){4}}
\put(47,11){\vector(1,0){4}}

\put(24,13){\vector(0,-1){14}}
\put(1,6){\vector(1,0){47}}

\put(22.5,-0.5){\scriptsize $i$}
\put(46,4.3){\scriptsize $j$}

\put(24.5,2.5){\scriptsize $2$}
\put(24.5,0.5){\scriptsize $3$}

\put(21.8,8.5){\scriptsize $-1$}
\put(21.8,10.5){\scriptsize $-2$}

\put(28.7,4.5){\scriptsize $2$}
\put(30.7,4.5){\scriptsize $3$}
\put(32.7,4.5){\scriptsize $4$}
\put(34.7,4.5){\scriptsize $5$}
\put(36.7,4.5){\scriptsize $6$}
\put(38.7,4.5){\scriptsize $7$}

\put(17.8,6.5){\scriptsize $-3$}
\put(15.8,6.5){\scriptsize $-4$}
\put(13.8,6.5){\scriptsize $-5$}
\put(11.8,6.5){\scriptsize $-6$}
\put(9.8,6.5){\scriptsize $-7$}
\put(7.8,6.5){\scriptsize $-8$}

\end{picture}
\end{center}
\setlength{\unitlength}{1pt}

\bigskip\bigskip
\noindent
Note that in the case $n=0$ the definition of the diagram formally
makes sense with $\la$ considered to be empty. Thus,
$\Ga(\emptyset)$ consists of two infinite rows of cells:

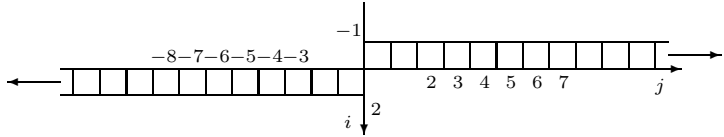
\setlength{\unitlength}{0.5em}
\begin{center}
\begin{picture}(48,10)

\put(1,4){\line(1,0){23}}
\put(24,8){\line(1,0){23}}

\multiput(2,4)(2,0){12}{\line(0,1){2}}
\multiput(24,6)(2,0){12}{\line(0,1){2}}

\put(1,5){\vector(-1,0){4}}
\put(47,7){\vector(1,0){4}}

\put(24,11){\vector(0,-1){10}}
\put(1,6){\vector(1,0){47}}

\put(22.5,1.5){\scriptsize $i$}
\put(46,4.3){\scriptsize $j$}

\put(24.5,2.5){\scriptsize $2$}

\put(21.8,8.5){\scriptsize $-1$}

\put(28.7,4.5){\scriptsize $2$}
\put(30.7,4.5){\scriptsize $3$}
\put(32.7,4.5){\scriptsize $4$}
\put(34.7,4.5){\scriptsize $5$}
\put(36.7,4.5){\scriptsize $6$}
\put(38.7,4.5){\scriptsize $7$}

\put(17.8,6.5){\scriptsize $-3$}
\put(15.8,6.5){\scriptsize $-4$}
\put(13.8,6.5){\scriptsize $-5$}
\put(11.8,6.5){\scriptsize $-6$}
\put(9.8,6.5){\scriptsize $-7$}
\put(7.8,6.5){\scriptsize $-8$}

\end{picture}
\end{center}
\setlength{\unitlength}{1pt}

\bigskip

By the {\it content\/} of a cell $\al=(i,j)$ with
coordinates $i$ and $j$ we shall mean the number $c(\al)=j-i$.
For any nonnegative integer $p$ we shall denote by $\Ga(\la)^{(p)}$
the image of $\Ga(\la)$ with respect to the shift operator
$(i,j)\mapsto (i-p,j)$. In other words, $\Ga(\la)^{(p)}$
is obtained from the diagram $\Ga(\la)$ by lifting each cell
$p$ units up.

We let $P_{1}(u),\dots,P_{n-m}(u)$ denote the Drinfeld polynomials
corresponding to the $\Y(\spa_{N-M})$-module $V(\lambda)^+_{\mu}$.

\bth\label{thm:dp}
For each $k=1,\dots,n-m$ the Drinfeld polynomial $P_k(u)$ is given by
\ben
P_k(u)=\prod_{\al}(u+c(\al)+1/2),
\een
where $\al$ runs over the cells of the intersection
$\Ga(\mu)\cap\Ga(\la)^{(k-1)}$.
\eth

\bex\label{ex:lamu}
Let $\la=(-2,-8,-10,-13)$ and $\mu=(-4,-7)$ so that $n=4$ and $m=2$.
(In a more standard notation, $\la$ and $\mu$ can be thought of as
the partitions $(13,10,8,2)$ and $(7,4)$, respectively).
The polynomial $P_1(u)$ is calculated
from the figure:

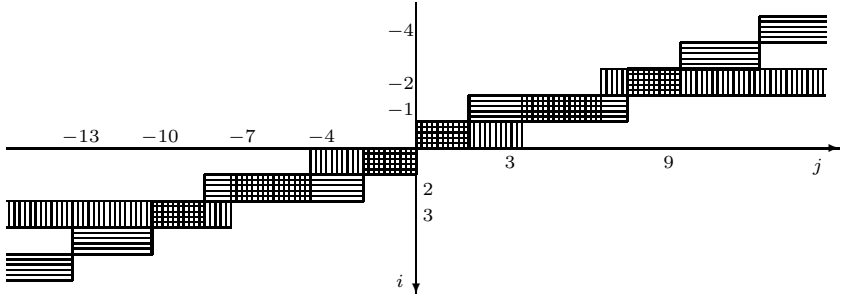
\setlength{\unitlength}{0.5em}
\begin{center}
\begin{picture}(64,22)

\put(1,4){\line(1,0){17}}
\put(1,6){\line(1,0){23}}
\put(18,8){\line(1,0){14}}
\put(24,8){\line(1,0){8}}
\put(24,10){\line(1,0){16}}
\put(32,10){\line(1,0){8}}
\put(32,12){\line(1,0){14}}
\put(40,12){\line(1,0){6}}
\put(40,14){\line(1,0){23}}
\put(46,16){\line(1,0){17}}

\thicklines

\put(1,0){\line(1,0){5}}
\put(1,2){\line(1,0){11}}
\put(6,2){\line(1,0){6}}
\put(6,4){\line(1,0){10}}
\put(12,6){\line(1,0){16}}
\put(16,8){\line(1,0){16}}
\put(28,10){\line(1,0){8}}
\put(32,12){\line(1,0){16}}
\put(36,14){\line(1,0){16}}
\put(48,16){\line(1,0){10}}
\put(52,18){\line(1,0){11}}
\put(58,20){\line(1,0){5}}

\put(6,0){\line(0,1){4}}
\put(12,2){\line(0,1){4}}
\put(16,4){\line(0,1){4}}
\put(28,6){\line(0,1){4}}
\put(32,8){\line(0,1){4}}
\put(36,10){\line(0,1){4}}
\put(48,12){\line(0,1){4}}
\put(52,14){\line(0,1){4}}
\put(58,16){\line(0,1){4}}

\thinlines

\put(32,21){\vector(0,-1){22}}
\put(1,10){\vector(1,0){63}}

\put(30.5,-0.5){\scriptsize $i$}
\put(62,8.3){\scriptsize $j$}

\multiput(1.2,4)(0.4,0){43}{\line(0,1){2}}
\multiput(18,6)(0.4,0){16}{\line(0,1){2}}
\multiput(24,8)(0.4,0){20}{\line(0,1){2}}
\multiput(32,10)(0.4,0){21}{\line(0,1){2}}
\multiput(40,12)(0.4,0){16}{\line(0,1){2}}
\multiput(46,14)(0.4,0){43}{\line(0,1){2}}

\multiput(1,0)(0,0.4){5}{\line(1,0){5}}
\multiput(6,2)(0,0.4){5}{\line(1,0){6}}
\multiput(12,4)(0,0.4){5}{\line(1,0){4}}
\multiput(16,6)(0,0.4){5}{\line(1,0){12}}
\multiput(28,8)(0,0.4){5}{\line(1,0){4}}
\multiput(32,10)(0,0.4){5}{\line(1,0){4}}
\multiput(36,12)(0,0.4){5}{\line(1,0){12}}
\multiput(48,14)(0,0.4){5}{\line(1,0){4}}
\multiput(52,16)(0,0.4){5}{\line(1,0){6}}
\multiput(58,18)(0,0.4){5}{\line(1,0){5}}

\put(32.5,6.5){\scriptsize $2$}
\put(32.5,4.5){\scriptsize $3$}

\put(29.8,12.5){\scriptsize $-1$}
\put(29.8,14.5){\scriptsize $-2$}
\put(29.8,18.5){\scriptsize $-4$}

\put(38.7,8.5){\scriptsize $3$}
\put(50.7,8.5){\scriptsize $9$}

\put(23.8,10.5){\scriptsize $-4$}
\put(17.8,10.5){\scriptsize $-7$}
\put(11.2,10.5){\scriptsize $-10$}
\put(5.2,10.5){\scriptsize $-13$}

\end{picture}
\end{center}
\setlength{\unitlength}{1pt}

\bigskip\bigskip
\noindent
The horizontal and vertical shadings indicate
the diagrams $\Ga(\la)$ and $\Ga(\mu)$,
respectively.
The cells belonging to the intersection $\Ga(\mu)\cap\Ga(\la)$
have the coordinates $(3,-10)$, $(3,-9)$, $(2,-7)$, $(2,-6)$, $(2,-5)$,
$(1,-2)$, $(1,-1)$, $(0,0)$, $(0,1)$, $(-1,4)$, $(-1,5)$, $(-1,6)$,
$(-2,8)$, $(-2,9)$. Hence,
\begin{multline}
P_1(u)=(u-25/2)(u-23/2)(u-17/2)(u-15/2)(u-13/2)(u-5/2)(u-3/2)\\
(u+1/2)(u+3/2)(u+11/2)(u+13/2)(u+15/2)(u+21/2)(u+23/2).
\non
\end{multline}
Note that the property $P_1(u)=P_1(-u+1)$ is implied by the central
symmetry of the set $\Ga(\mu)\cap\Ga(\la)$.

The polynomial $P_2(u)$ is calculated
from the figure:

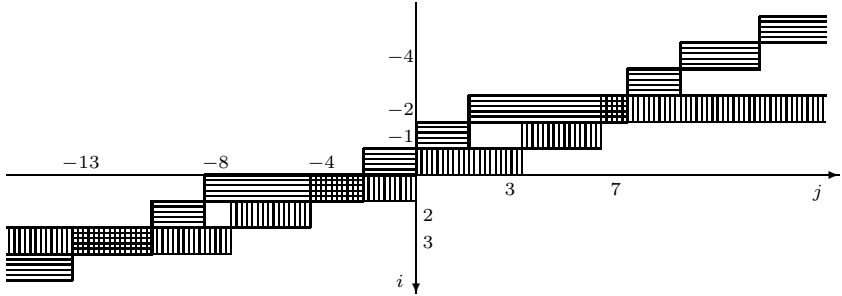
\setlength{\unitlength}{0.5em}
\begin{center}
\begin{picture}(64,24)

\put(1,4){\line(1,0){17}}
\put(1,6){\line(1,0){23}}
\put(18,8){\line(1,0){14}}
\put(24,8){\line(1,0){8}}
\put(24,10){\line(1,0){16}}
\put(32,10){\line(1,0){8}}
\put(32,12){\line(1,0){14}}
\put(40,12){\line(1,0){6}}
\put(40,14){\line(1,0){23}}
\put(46,16){\line(1,0){17}}

\thicklines

\put(1,2){\line(1,0){5}}
\put(1,4){\line(1,0){11}}
\put(6,6){\line(1,0){6}}
\put(6,6){\line(1,0){10}}
\put(12,8){\line(1,0){16}}
\put(16,10){\line(1,0){16}}
\put(28,12){\line(1,0){8}}
\put(32,14){\line(1,0){16}}
\put(36,16){\line(1,0){16}}
\put(48,18){\line(1,0){10}}
\put(52,20){\line(1,0){11}}
\put(58,22){\line(1,0){5}}

\put(6,2){\line(0,1){4}}
\put(12,4){\line(0,1){4}}
\put(16,6){\line(0,1){4}}
\put(28,8){\line(0,1){4}}
\put(32,10){\line(0,1){4}}
\put(36,12){\line(0,1){4}}
\put(48,14){\line(0,1){4}}
\put(52,16){\line(0,1){4}}
\put(58,18){\line(0,1){4}}

\thinlines

\put(32,23){\vector(0,-1){22}}
\put(1,10){\vector(1,0){63}}

\put(30.5,1.5){\scriptsize $i$}
\put(62,8.3){\scriptsize $j$}

\multiput(1.2,4)(0.4,0){43}{\line(0,1){2}}
\multiput(18,6)(0.4,0){16}{\line(0,1){2}}
\multiput(24,8)(0.4,0){20}{\line(0,1){2}}
\multiput(32,10)(0.4,0){21}{\line(0,1){2}}
\multiput(40,12)(0.4,0){16}{\line(0,1){2}}
\multiput(46,14)(0.4,0){43}{\line(0,1){2}}

\multiput(1,2)(0,0.4){5}{\line(1,0){5}}
\multiput(6,4)(0,0.4){5}{\line(1,0){6}}
\multiput(12,6)(0,0.4){5}{\line(1,0){4}}
\multiput(16,8)(0,0.4){5}{\line(1,0){12}}
\multiput(28,10)(0,0.4){5}{\line(1,0){4}}
\multiput(32,12)(0,0.4){5}{\line(1,0){4}}
\multiput(36,14)(0,0.4){5}{\line(1,0){12}}
\multiput(48,16)(0,0.4){5}{\line(1,0){4}}
\multiput(52,18)(0,0.4){5}{\line(1,0){6}}
\multiput(58,20)(0,0.4){5}{\line(1,0){5}}

\put(32.5,6.5){\scriptsize $2$}
\put(32.5,4.5){\scriptsize $3$}

\put(29.8,18.5){\scriptsize $-4$}
\put(29.8,14.5){\scriptsize $-2$}
\put(29.8,12.5){\scriptsize $-1$}

\put(38.7,8.5){\scriptsize $3$}
\put(46.7,8.5){\scriptsize $7$}

\put(23.8,10.5){\scriptsize $-4$}
\put(15.8,10.5){\scriptsize $-8$}
\put(5.2,10.5){\scriptsize $-13$}

\end{picture}
\end{center}
\setlength{\unitlength}{1pt}

\noindent
The cells which belong to the intersection $\Ga(\mu)\cap\Ga(\la)^{(1)}$
have the coordinates $(3,-13)$, $(3,-12)$, $(3,-11)$,
$(1,-4)$, $(1,-3)$, $(-2,7)$. Hence,
\ben
P_2(u)=(u-31/2)(u-29/2)(u-27/2)(u-9/2)(u-7/2)(u+19/2).
\een
\eex

\bex\label{ex:eval}
In the case $m=0$ the vector space $V(\la)^+_{\mu}$
can be identified with $V(\la)$ and the corresponding
$\Y(\spa_{2n})$-module coincides with
the evaluation module defined by the evaluation homomorphism \eqref{evaltw}.
Applying Theorem~\ref{thm:dp} to the diagrams $\Ga(\la)$ and $\Ga(\emptyset)$
we obtain
\begin{multline}
P_1(u)=(u+\la_1-1/2)(u+\la_1+1/2)\cdots (u-3/2)\\
{}\times(u+1/2)(u+3/2)\cdots (u-\la_1-1/2)
\non
\end{multline}
and
\ben
P_k(u)=(u+\la_{k}-1/2)(u+\la_{k}+1/2)\cdots (u+\la_{k-1}-3/2), \qquad
k=2,\dots,n.
\een
On the other hand,
the highest weight of this module is given by \eqref{laimui}.
Due to \eqref{dpk}
and \eqref{dpfirst} this
obviously agrees with the above calculation of the $P_k(u)$.
\eex

\begin{proof}[Proof of Theorem~\ref{thm:dp}]
We shall derive the statement from Theorem~\ref{thm:hwmu} and
the definition of the Drinfeld polynomials \eqref{dpk}
and \eqref{dpfirst}. In order to calculate $P_{1}(u)$ observe that
the component $\mu_{m+1}(u)$ of the highest weight can also be found
from the formula
\ben
\vp_{m+1,m+1}(u)\ts \zeta^{}_{\La_0}=\mu_{m+1}(u)\ts \zeta^{}_{\La_0};
\een
see \eqref{phiab}. Hence, applying \eqref{alfkk} with $k=m+1$
and using \eqref{maxmid} and \eqref{minmid} we get
\begin{multline}
\mu_{m+1}(u)=
\prod_{i=1}^{m+1}\frac{u-\middle\{\la_{i-1},\mu_{i-1},\mu_{i}\}+i-1/2}
{u+i-1/2}\\
{}\times
\prod_{i=1}^{m+1}\frac{u+\middle\{\la_{i},\mu_{i-1},\mu_{i}\}-i+1/2}
{u-i+1/2}.
\non
\end{multline}
Therefore, by \eqref{dpfirst},
\begin{multline}
\frac{P_{1}(u+1)}{P_{1}(u)}=
\prod_{i=1}^{m+1}\frac{u+\middle\{\la_{i-1},\mu_{i-1},\mu_{i}\}-i+1/2}
{u+\middle\{\la_{i},\mu_{i-1},\mu_{i}\}-i+1/2}\\
{}\times
\prod_{i=1}^{m+1}\frac{u-\middle\{\la_{i},\mu_{i-1},\mu_{i}\}+i-1/2}
{u-\middle\{\la_{i-1},\mu_{i-1},\mu_{i}\}+i-1/2}.
\non
\end{multline}
This gives
$P_{1}(u)=Q(u)\ts Q(-u+1)\ts (-1)^{\deg Q}$, where
\beql{qu}
Q(u)=\prod_{i=1}^{m+1}(u+\be_i)(u+\be_i+1)\cdots (u+\al_i-1)
\eeq
with
\ben
\al_i=\middle\{\la_{i-1},\mu_{i-1},\mu_{i}\}-i+1/2\Fand
\be_i=\middle\{\la_{i},\mu_{i-1},\mu_{i}\}-i+1/2.
\een
Thus, the factor corresponding to the index $i$
in the product in \eqref{qu} can be interpreted as the product
$\prod_{\al}(u+c(\al)+1/2)$, where $\al$ runs over
the cells of the intersection of the diagrams $\Ga(\la)\cap\Ga(\mu)$
whose first coordinate is $i$. Taking into account the central symmetry
of $\Ga(\la)\cap\Ga(\mu)$, we derive the desired formula for $P_{1}(u)$.

Now using \eqref{dpk} and replacing the index $k$ by
$r=k-m$ in Theorem~\ref{thm:hwmu}
we derive the following
expression for the Drinfeld polynomial $P_{r+1}(u)$ with
$1\leqslant r< n-m$:
\beql{pkpk}
\bal
\frac{P_{r+1}(u+1)}{P_{r+1}(u)}&=
\prod_{\overset{\scl{i=1}}{\la_i<\mu_{i+r}}}^{r+m-2}
\frac{u-\max\{\la_{i+1},\mu_{i+r+1}\}+r+i+1/2}
{u-\max\{\la_{i},\mu_{i+r+1}\}+r+i+1/2}\\
{}&\times
\prod_{\overset{\scl{i=1}}{\la_i<\mu_{i+r-1}\leqslant\la_{i-1}}}^{r+m-1}
\frac{u-\max\{\la_{i},\mu_{i+r}\}+r+i-1/2}{u-\mu_{i+r-1}+r+i-1/2}\\
{}&\times
\prod_{\overset{\scl{i=2}}{\la_i>\mu_{i-r}}}^{r+m+1}
\frac{u+\min\{\la_{i-1},\mu_{i-r-1}\}+r-i+1/2}
{u+\min\{\la_{i},\mu_{i-r-1}\}+r-i+1/2}\\
{}&\times
\prod_{\overset{\scl{i=1}}{\la_{i+1}\leqslant\mu_{i-r+1}<\la_i}}^{r+m-1}
\frac{u+\min\{\la_{i},\mu_{i-r}\}+r-i-1/2}
{u+\mu_{i-r+1}+r-i+1/2}.
\eal
\eeq
Note that the expression has the form
\ben
\frac{P_{r+1}(u+1)}{P_{r+1}(u)}=\prod_j\frac{u+\al_j}{u+\be_j}
\een
for some parameters $\al_j$ and $\be_j$ with $\al_j-\be_j\in\ZZ_+$.
This implies that $P_{r+1}(u)$ is given by
\beql{pkbeal}
P_{r+1}(u)=\prod_j(u+\be_j)(u+\be_j+1)\cdots (u+\al_j-1).
\eeq
It is straightforward to verify that this product
coincides with $\prod_{\al}(u+c(\al)+1/2)$, where $\al$ runs over
the cells of the intersection $\Ga(\mu)\cap\Ga(\la)^{(r-1)}$.
Indeed, changing the product index $i$
respectively by $i-r$ and by $i-r+1$ in the third and fourth
products in \eqref{pkpk} we bring these products to the form
\begin{multline}
\label{prred}
\prod_{\overset{\scl{i=1}}{\la_{i+r}>\mu_{i}}}^{m+1}
\frac{u+\min\{\la_{i+r-1},\mu_{i-1}\}-i+1/2}
{u+\min\{\la_{i+r-1},\mu_{i-1}\}-i+1/2}\\
{}\times
\prod_{\overset{\scl{i=1}}{\la_{i+r}\leqslant\mu_{i}<\la_{i+r-1}}}^{m}
\frac{u+\min\{\la_{i+r-1},\mu_{i-1}\}-i+1/2}
{u+\mu_{i}-i+1/2}.
\end{multline}
Therefore, the factor corresponding to the index $i$
contributes into \eqref{pkbeal} the product
$\prod_{\al}(u+c(\al)+1/2)$, where $\al$ runs over
the cells of the intersection $\Ga(\mu)\cap\Ga(\la)^{(r-1)}$
whose first coordinate is $i$. Note that writing
\ben
-\max\{\la_{i},\mu_{i+r+1}\}=\min\{\la_{-i},\mu_{-i-r-1}\}
\een
we can similarly bring the first and second
products in \eqref{pkpk} to the form \eqref{prred},
where the product is taken over negative indices.
This
contributes the product $\prod_{\al}(u+c(\al)+1/2)$ into \eqref{pkbeal},
where $\al$ runs over
the cells of the intersection $\Ga(\mu)\cap\Ga(\la)^{(r-1)}$
whose first coordinate is non-positive.
\epf

\section{Centralizer construction}
\label{sec:cc}
\setcounter{equation}{0}

We now return to our general situation so that $\g_N$ denotes either the
orthogonal or symplectic Lie algebra; see \eqref{oos}.
We start by recalling the construction of the {\it Olshanski algebra\/} $\Ar_M$;
see \cite{o:ty, mo:cc}.
Fix a nonnegative integer $M$ such that $N-M$ is even.
So, if $N=2n$ or $N=2n+1$ then
$M=2m$ or $M=2m+1$,
respectively, for some $m\leqslant n$.
Denote by $\ggot'_{N-M}$
the subalgebra of $\g_N$ spanned by
the elements $F_{ij}$ subject to the condition $m+1\leq |i|,|j|\leq n$.
Let $\Ar_M(N)$ denote the centralizer
of $\ggot'_{N-M}$ in the universal enveloping algebra $\U(\g_N)$.
Let $\U(\g_N)^0$ denote the centralizer of $F_{nn}$
in $\U(\g_N)$ and let $\Ir(N)$ be the left ideal in $\U(\g_N)$ generated by the
elements $F_{in}$, $i=-n,\dots,n$. Then $\Ir(N)^0=\Ir(N)\cap \U(\g_N)^0$
is a two-sided ideal in $\U(\g_N)^0$
which coincides with the intersection $\Jr(N)\cap \U(\g_N)^0$,
where $\Jr(N)$ is the right ideal in $\U(\g_N)$ generated by the
elements $F_{ni}$, $i=-n,\dots,n$.
One has a vector space
decomposition
\ben
\U(\g_N)^0=\Ir(N)^0\oplus \U(\g_{N-2}).
\een
Therefore the projection of $\U(\g_N)^0$ onto $\U(\g_{N-2})$
with the kernel $\Ir(N)^0$ is an algebra homomorphism. Its
restriction to the subalgebra $\Ar_M(N)$ defines
a filtration-preserving homomorphism
\beql{pinaa}
\pi_N: \Ar_M(N)\to \Ar_M(N-2)
\eeq
so that one can define the algebra $\Ar_M$ as
the projective limit with respect to this sequence of homomorphisms
in the category of filtered algebras; see \cite{mo:cc} for more details.

Taking the composition of the homomorphism $\Y(\g_M)\to\Y(\g_N)$
defined in Proposition~\ref{prop:dualhomtw} and the evaluation
homomorphism \eqref{evaltw} we obtain another homomorphism
$\psi_N:\Y(\g_M)\to\U(\g_{N})$
which takes $s_{ij}(u)$ to the series
\beql{alfij}
\al_{m-n}(u)\cdot
{f}^{-n\ts\cdots\ts -m-1,\ts i,\ts m+1\cdots \ts n}_{-n\ts\cdots\ts -m-1,
\ts j,\ts m+1\cdots \ts n}\ts(u+n-m),\qquad -m\leqslant i,j\leqslant m,
\eeq
where we have used the notation of Section~\ref{sec:sr} for the images
of the Sklyanin minors with respect to \eqref{evaltw}. Observe now
that by Corollary~\ref{cor:centsm} the image of $\psi_N$ is contained
in the centralizer $\Ar_M(N)$, so that
$\psi_N:\Y(\g_M)\to\Ar_M(N)$.

\bpr\label{prop:comd}
The sequence of homomorphisms
$(\psi_N|\ N=M+2k, k=0,1,\dots)$ defines a homomorphism
\ben
\psi: \Y(\g_M)\to\Ar_M.
\een
\epr

\bpf
We have to verify that the homomorphisms $\psi_N$
are compatible with the sequence of homomorphisms \eqref{pinaa}, that is, the
following diagram is commutative:
\ben
\begin{CD}
\Y(\g_M) @= \Y(\g_M) @= \cdots @= \Y(\g_M)@=\cdots\\
@V \psi_M VV     @V\psi_{M+2}VV @.
@V\psi_{N}VV\\
\Ar_M(M) @<<\pi_{M+2}< \Ar_{M}(M+2)@<<<\cdots
@<<\pi_{N}< \Ar_M(N)@<<<\cdots.
\end{CD}
\een
Let us calculate the image of the series $\psi_N(s_{ij}(u))$ under
the homomorphism $\pi_N$. Applying \eqref{sscheck} we obtain
\begin{multline}
{f}^{-n\ts\cdots\ts -m-1,\ts i,\ts m+1\cdots \ts n}_{-n\ts\cdots\ts -m-1,
\ts j,\ts m+1\cdots \ts n}\ts(u+n-m)\\
{}=\sum_c
{\check f}^{-n\ts\cdots\ts -m-1,\ts i,\ts m+1\cdots \ts n}_{-n\ts\cdots\ts -m-1,
\ts j,\ts m+1\cdots \ts n-1,c}\ts(u+n-m)\ts f_{cn}(u-n+m),
\non
\end{multline}
where by ${\check f\ts}^{a_1\cdots\ts a_k}_{b_1\cdots\ts b_{k-1},c}(u)$
we denote the image of the auxiliary minor
${\check s\ts}^{a_1\cdots\ts a_k}_{b_1\cdots\ts b_{k-1},c}(u)$
under the evaluation homomorphism \eqref{evaltw}.
Now observe that $f_{cn}(u-n+m)$ belongs to the left ideal $\Ir(N)$
unless $c=n$. In this case $f_{nn}(u-n+m)\equiv 1\mod\Ir(N)$.
Furthermore, by \eqref{auzerosf},
\begin{multline}
{\check f}^{-n\ts\cdots\ts -m-1,\ts i,\ts m+1\cdots \ts n}_{-n\ts\cdots\ts -m-1,
\ts j,\ts m+1\cdots \ts n}\ts(v)\\
=\frac{2v+1}{2v\pm1}\ts\sum_{i=1}^{2n-2m}(-1)^{i-1}\ts
f^t_{a_i,-n}(-v)
\ts {f\ts}^{\ a_1\cdots\ts\wh a_i\ts\cdots \ts a_{2n-2m}}_{-n+1\ts\cdots\ts -m-1,
\ts j,\ts m+1\cdots \ts n-1}(v-1),
\non
\end{multline}
where we have set $(a_1,\dots,a_{2n-2m})=
(-n\ts\cdots\ts -m-1,\ts i,\ts m+1\cdots \ts n-1)$ and $v=u+n-m$.
Now $f^t_{a_i,-n}(-v)=\theta_{-n,a_i}\ts f_{n,-a_i}(-v)$
belongs to the right ideal $\Jr(N)$ unless $a_i=-n$, that is, $i=1$.
In this case $f_{nn}(-v)\equiv 1\mod\Jr(N)$. Note that
\ben
\al_{m-n}(u)\ts \frac{2v+1}{2v\pm1}=\al_{m-n+1}(u).
\een
Hence, $\pi_N$ takes $\psi_N(s_{ij}(u))$ to the series
\ben
\al_{m-n+1}(u)\cdot
{f}^{-n+1\ts\cdots\ts -m-1,\ts i,\ts m+1\cdots \ts n-1}_{-n+1\ts\cdots\ts -m-1,
\ts j,\ts m+1\cdots \ts n-1}\ts(u+n-m-1)
\een
which coincides with $\psi_{N-2}(s_{ij}(u))$. Thus, the sequence
of coefficients at each power of $u^{-1}$ in
$\psi_N(s_{ij}(u))$ with $N=M+2k$, $k=0,1,\dots$ defines an element of $\Ar_M$.
\epf

Corollary~\ref{cor:centsm} implies that all the coefficients of the series
\beql{fminmn}
\al_{m-n}(u)\cdot
{f}^{-n\ts\cdots\ts -m-1,\ts m+1\cdots \ts n}_{-n\ts\cdots\ts -m-1,
\ts m+1\cdots \ts n}\ts(u+n-m)
=1+c^{(n)}_1u^{-1}+c^{(n)}_2u^{-2}+\cdots
\eeq
belong to the center of the universal enveloping algebra $\U(\ggot'_{N-M})$.
Hence, for any $i$ the coefficient $c^{(n)}_i$ is an element
of the centralizer $\Ar_M(N)$. The argument of the proof of
Proposition~\ref{prop:comd} shows that the image of the series
\eqref{fminmn} under the homomorphism $\pi_N$ is
\ben
\al_{m-n+1}(u)\cdot
{f}^{-n+1\ts\cdots\ts -m-1,\ts m+1\cdots \ts n-1}_{-n+1\ts\cdots\ts -m-1,
\ts m+1\cdots \ts n-1}\ts(u+n-m-1).
\een
Therefore, for each $i$ the sequence $c_i=(c^{(n)}_i\ |\ n\geqslant m+1)$
determines an element of the projective limit algebra $\Ar_M$.
Denote by $\Cr_M$ the subalgebra of $\Ar_M$ generated by all
$c_i$ with $i\geq 1$. The subalgebra $\Cr_M$ is
studied in detail in \cite[Section~3]{mo:cc} where
it was identified with
the algebra of {\it virtual Laplace operators\/}.
Up to an obvious
change of notation,
the series \eqref{fminmn}
coincides with the Sklyanin minor
$f_{\Bc}(u)$ introduced in the proof of Theorem~\ref{thm:hwmu}.
So its Harish-Chandra image can be found from \eqref{fbhc} which shows that
the elements $c_{2i}$ with even indices
are algebraically independent and generate the algebra $\Cr_M$;
cf. \cite[Section~3]{mo:cc}.

\bth\label{thm:cctw}
The homomorphism
$
\psi: \Y(\g_M)\hookrightarrow\Ar_M
$
is injective.
Moreover, one has
an isomorphism
\ben
\Ar_M=\Cr_M\ot \Y(\g_M),
\een
where $\Y(\g_M)$ is identified with its image under the embedding $\psi$.
\eth

\bpf
Our argument is similar to the proof of \cite[Theorem~4.17]{mo:cc}.
Consider the canonical filtration of the universal
enveloping algebra $\U(\g_{N})$. The corresponding
graded algebra $\gr \U(\g_{N})$ is isomorphic to
the symmetric algebra
$\Sr(\g_{N})$ of the space
$\g_N$. Elements of $\Sr(\g_{N})$ can be naturally identified
with polynomials in matrix elements of an $N\times N$ matrix
$x=(x_{ij})_{i,j=-n}^n$ such that $x^t=-x$.
Denote by
$\Prm_M(N)$ the subalgebra of the elements
of $\Sr(\g_{N})$ which are invariant under the adjoint action of
the Lie algebra $\ggot'_{N-M}$.
The algebra $\Ar_M$ possesses a natural filtration induced by
the canonical filtrations on the centralizers $\Ar_M(N)$.
The corresponding graded algebra $\gr \Ar_M$ is naturally isomorphic to
the projective limit $\Prm_M$ of the commutative algebras $\Prm_M(N)$
with respect to homomorphisms $\Prm_M(N)\to \Prm_M(N-2)$
analogous to \eqref{pinaa}; see \cite[Section~4]{mo:cc}.
The images
in $\Prm_M(N)$ of the coefficients the series \eqref{alfij} and
\eqref{fminmn} can be found from the explicit formulas for the Sklyanin
minors; see Proposition~\ref{prop:minsf}.
Indeed, apply the proposition to the Sklyanin minor
${s\ts}^{-n\ts\cdots\ts -m-1,\ts m+1\cdots \ts n}_{-n\ts\cdots\ts -m-1,
\ts m+1\cdots \ts n}\ts(u+n-m)$
then replace each series $s_{ij}(u)$ by its image $f_{ij}(u)$
under the evaluation homomorphism \eqref{evaltw}.
Observe that the image of $s^t_{ij}(-u)$ coincides with $f_{ij}(u\mp1)$.
Since we are only interested in the highest degree component of
the coefficient at each power of $u^{-1}$, we may replace
each expression of type $f_{ij}(u+c)$ by $\delta_{ij}+F_{ij}u^{-1}$.
Hence,
denoting the elements of the set $\Ac=\{-n,\dots,-m-1,m+1,\dots,n\}$
by $a_1,\dots, a_{N-M}$ we can write,
modulo lower degree terms at each power of $u$,
\begin{multline}
{f\ts}^{-a_1\ts\cdots \ts -a_{N-M}}_{a_1\ts\cdots\ts a_{N-M}}(u+n-m)
\equiv\al_{n-m}(u+n-m)\\
{}\times\sum_{p\in\Sym_M}
\sgn p\tss p'\cdot (1+Fu^{-1})_{-a_{p(1)},a_{p'(1)}}\cdots
(1+Fu^{-1})_{-a_{p(N-M)},a_{p'(N-M)}},
\non
\end{multline}
where $F$ denotes the $(N-M)\times (N-M)$ matrix whose rows and columns
are enumerated by the elements of the set $\Ac$ and whose $ij$-th entry
is $F_{ij}$. Since these matrix elements commute modulo lower degree
terms, taking into account the relation
\ben
\al_{m-n}(u)\ts\al_{n-m}(u+n-m)=1,
\een
we can conclude from Lemma~\ref{lem:probi} that the image of the series
\eqref{fminmn} in $\Prm_M(N)$ is the determinant
$
\det (1+x u^{-1})^{}_{\Ac\Ac}.
$
The same argument shows that the image of the series
\eqref{alfij} in $\Prm_M(N)$ is
$
\det (1+x u^{-1})^{}_{\Ac_i\Ac_j},
$
where
\ben
\Ac_i=\{-n,\dots,-m-1,i,m+1,\dots,n\}.
\een
Our next step is to show that every element $\phi$
of the algebra $\Prm_M(N)$ such that $\deg\phi<n-m$
can be represented as a polynomial in the coefficients of
the series $\det (1+x u^{-1})^{}_{\Ac\Ac}$ and
$\det (1+x u^{-1})^{}_{\Ac_i\Ac_j}$. However, it
was proved in \cite[Section~4.9]{mo:cc} that $\phi$
can be presented as a polynomial in the elements
$\tr (x_{\Ac\Ac})^k$ and $\La^{(k)}_{ij}$
with $-m\leqslant i,j\leqslant m$ and $k\geqslant 1$,
where $\La^{(k)}_{ij}=\sum x_{ic_1}x_{c_1c_2}\cdots x_{c_{k-1}j}$,
summed over the indices $c_r\in\Ac$. On the other hand,
each element $\tr (x_{\Ac\Ac})^k$ is a polynomial in the
coefficients of the series $\det (1+x u^{-1})^{}_{\Ac\Ac}$.
This follows from the fact that the coefficients
of the characteristic polynomial $\det (u+x)^{}_{\Ac\Ac}$
form a complete set of invariants of the matrix $x_{\Ac\Ac}$.
(Explicit expressions for the elements $\tr (x_{\Ac\Ac})^k$
in terms of the coefficients of the characteristic polynomial
can be derived e.g. from the Liouville formula;
cf. \cite[Remark~5.8]{mno:yc}). The claim now follows from
the well known identity
\beql{quasid}
\det (1+x u^{-1})^{}_{\Ac_i\Ac_j}=\det (1+x u^{-1})^{}_{\Ac\Ac}
\cdot
\Big(\delta_{ij}+\sum_{k=1}^{\infty}(-1)^{k-1}\Lambda_{ij}^{(k)} u^{-k}\Big);
\eeq
see e.g. \cite{gr:dm,gr:tn}. Indeed, the identity implies that
the elements $\Lambda_{ij}^{(k)}$ are polynomials
in the coefficients of the series $\det (1+x u^{-1})^{}_{\Ac_i\Ac_j}$
and $\det (1+x u^{-1})^{}_{\Ac\Ac}$. This allows us to conclude that
the algebra $\Ar_M$ is generated by the subalgebra $\Cr_M$ and the image
of the homomorphism $\psi$.

Observe that since the matrix $x$ satisfies $x^t=-x$ we have the
relations
\ben
\det (1+x u^{-1})^{}_{\Ac\Ac}=\det (1-x u^{-1})^{}_{\Ac\Ac}
\een
and
\ben
\Lambda_{ij}^{(k)}=(-1)^k\ts\theta_{ij}\ts\Lambda_{-j,-i}^{(k)}.
\een
Hence, we can write
\ben
\det (1+x u^{-1})^{}_{\Ac\Ac}=1+\sum_{r=1}^{\infty}\La^{(2r)}\ts u^{-2r}
\een
for some polynomials $\La^{(2r)}$ in the matrix elements of $x$.
For the elements $\Lambda_{ij}^{(k)}$ we shall impose
the following restrictions on $i,j,k$:
\ben
i+j<0\quad\text{for}\ \  k \ \ \text{odd,}\qquad
i+j\leq 0\quad\text{for}\ \  k \ \ \text{even}
\een
in the orthogonal case, and
\ben
i+j<0\quad \text{for}\ \  k \ \ \text{even,}\qquad
i+j\leq 0\quad\text{for}\ \  k \ \ \text{odd}
\een
in the symplectic case. Fix a positive integer $K$ and assume that the index
$k$ satisfies $1\leqslant k\leqslant K$.
It follows from the argument of
\cite[Section~4.10]{mo:cc} that there exists
a large enough value of $N$ such that the polynomials
$\La^{(k)}$ for even $k$, and $\Lambda_{ij}^{(k)}$ with the above
restrictions on $i,j,k$ are algebraically independent.
Due to the identity \eqref{quasid}, the same statement will hold
if each polynomial $\Lambda_{ij}^{(k)}$ is replaced by
the coefficient at $u^{-k}$ of the series
$\det (1+x u^{-1})^{}_{\Ac_i\Ac_j}$.
This proves the
injectivity of $\psi$ and the
tensor product decomposition for the algebra $\Ar_M$.
\epf

\end{document}